\documentclass{article}
\usepackage[utf8]{inputenc}
\usepackage[table]{xcolor}
\usepackage{graphicx}
\usepackage{color}
\usepackage{float}
\usepackage{todonotes}
\usepackage{amsfonts}
\usepackage{amsthm}
\usepackage{amsmath}
\usepackage{caption}
\usepackage{amssymb}
\usepackage{stmaryrd}
\usepackage{mathtools}
\usepackage[top=2.5cm, bottom=2.5cm, left=3cm, right=3cm]{geometry}
\usepackage[english]{babel}
\usepackage[T1]{fontenc}
\usepackage{hyperref}
\usepackage[ruled,vlined]{algorithm2e}
\usepackage{algorithmicx}
\usepackage{lineno}
\usepackage{float}
\usepackage{csquotes}
\usepackage{colortbl}
\usepackage{comment}

\newcommand{\R}{\mathbb{R}}

\newcommand{\E}{\mathbb{E}}

\newcommand{\di}{\mathrm{d}}

\def\Pc{{\cal P}}
\def\Lc{{\cal L}}
\def\red#1{{\color{red}#1}}

\def\Yc{{\cal Y}}
\def\Zc{{\cal Z}}
\def\Dc{{\cal D}}

\def\Fc{{\cal F}}
\def\Hc{{\cal H}}
\def\Kc{{\cal K}}
\def\Uc{{\cal U}}
\def\Tc{{\cal T}}
\def\Ac{{\cal A}}
\def\Lc{{\cal L}}

\def\Nc{{\cal N}}

\def \trans{^{\scriptscriptstyle{\intercal}}}
\DeclareMathOperator*{\argmin}{arg\,min}

\def \trans{^{\scriptscriptstyle{\intercal}}}

\def \E{\mathbb{E}}
\def \F{\mathbb{F}}

\def \P{\mathbb{P}}

\def\bop{{\boldsymbol p}}

\def\argmax{\mathop{\rm argmax}}
\def\argmin{\mathop{\rm argmin}}
\def\argmin_#1{\underset{#1}{\mathrm{argmin\, }}}

\def \d{\mathrm{d}} 

\def \Sum{\displaystyle\sum}

\def \Frac{\displaystyle\frac}

\def \b1{\bf{1}}

\def \mfa{\mathfrak{a}} 
\def \mra{\mathrm{a}} 
\def \mrJ{\mathrm{J}} 
 
\def \d{\mathrm{d}} 
 
\def \mrp{\mathrm{p}}

\def\beqs{\begin{eqnarray*}}
\def\enqs{\end{eqnarray*}}
\def\beq{\begin{eqnarray}}
\def\enq{\end{eqnarray}}

\def\argmin_#1{\underset{#1}{\mathrm{argmin\, }}}
\def\argmax_#1{\underset{#1}{\mathrm{argmax\, }}}

\def\bl#1{{\color{blue}#1}}

\numberwithin{equation}{section} 

\mathtoolsset{showonlyrefs}

\usepackage[backend=bibtex, maxbibnames=5, style=numeric]{biblatex} 
\begin{document}

\title{Mean-field neural networks-based  algorithms for McKean-Vlasov control problems
\thanks{This work is supported by  FiME, Laboratoire de Finance des March\'es de l'Energie, and the ``Finance and Sustainable Development'' EDF - CACIB Chair.}
}

\author{Huy\^en \sc{Pham}\footnote{LPSM,  Universit\'e Paris Cité,  \& FiME \sf \href{mailto:pham at lpsm.paris}{pham at lpsm.paris}} 
\and
Xavier \sc{Warin}
\footnote{EDF R\&D \& FiME \sf \href{mailto:xavier.warin at edf.fr}{xavier.warin at edf.fr}}}

\maketitle

\begin{abstract}
This paper is devoted to the numerical resolution of McKean-Vlasov control problems via the class of mean-field neural networks introduced in our companion paper \cite{phawar22a} in order to learn the solution on the Wasserstein space. 
We propose several algorithms either based on dynamic programming with control learning by policy or value iteration, or  backward SDE from stochastic maximum principle with global or local loss functions. Extensive numerical results on diffe\-rent examples are presented to illustrate the accuracy of 
each of our eight algorithms. We discuss and compare the pros and cons of all the tested methods.  
\end{abstract}


\vspace{5mm}


\section{Introduction}

This paper is concerned with the numerical resolution of  McKean-Vlasov (MKV)  control, also called mean-field control (MFC)  problems over finite horizon. 
The dynamics of the controlled state process $X$ $=$ $(X_t)_t$ valued in $\R^d$ is driven by the mean-field SDE (stochastic differential equation):
\beqs
\d X_t &=& b(X_t,\P_{X_t},\alpha_t) \d t + \sigma(X_t,\P_{X_t},\alpha_t) \d W_t,  \quad 0 \leq t\leq T, \quad X_0 \sim \mu_0, 
\enqs
where $W$ is a  $d$-dimensional Brownian motion on a filtered probability space $(\Omega,\Fc,\F=(\Fc_t)_t,\P)$, the initial distribution $\mu_0$ of $X_0$ lies in $\Pc_2(\R^d)$, the Wasserstein space of square-integrable probability measures, 
$\alpha$ $\in$ $\Ac$ is a control process, i.e,  an $\F$-progressively measurable process valued in $A$ $\subset$ $\R^m$, and $\P_{X_t}$ denotes the law of $X_t$, valued on 
$\Pc_2(\R^d)$, under standard assumptions on the coefficients $b$, $\sigma$ defined on $\R^d\times\Pc_2(\R^d)\times A$, and valued respectively in $\R^d$ and $\R^{d\times d}$. The objective is to minimize over controls $\alpha$ $\in$ $\Ac$, a cost functional of the form
\begin{align} \label{MKVcontrol}
J(\alpha) &=\;  \E \Big[ \int_0^T f(X_t,\P_{X_t},\alpha_t) \d t + g(X_T,\P_{X_T}) \Big], \quad \rightarrow \quad v(\mu_0) \; = \; \inf_{\alpha\in \Ac} J(\alpha),  
\end{align}
where $f$ is a running cost function on $\R^d\times\Pc_2(\R^d)\times A$, and $g$ is a terminal cost  function on $\R^d\times\Pc_2(\R^d)$.  

The theory and applications of mean-field control problems that study  models of  large population of interacting agents controlled by a social planner,  have  generated a vast literature in the last decade, and we refer to the monographs \cite{benetal13}, \cite{cardel19}, \cite{cardel2} for a comprehensive treatment of this topic. As analytical solutions to MFC are rarely available,  it is crucial to design efficient numerical schemes for solving such problem,  and the main challenging issue is the infinite dimensional feature of MFC coming from the distribution law state variable.

Following the tremendous impact of  machine learning methods for solving high-dimensional partial differential equations (PDEs) and control problems, see e.g. the survey papers \cite{bechutjenkuc20}, \cite{GPW21}, and the link to the website \href{http://deeppde.org/intro/}{deeppde.org},  some recent works have proposed deep learning schemes for MFC, based on neural network approximations of the feedback control and/or the value function solution to the Hamilton-Jacobi-Bellman equation or backward stochastic differential equations (BSDEs). In these articles, the authors consider either 
 approximate  feedback controls by standard feedforward neural networks with input the time and the state variable $X_t$ in  $\R^d$ by viewing the law of $X_t$ as a deterministic function of time (see \cite{pfe17}, \cite{carlau19},   \cite{fouquezhang19}, \cite{germicwar19}, \cite{Ruthotto9183}, \cite{reistozha21}),   or consider a particle approximation of the MFC for reducing the problem to a finite-dimensional problem that is numerically solved by means of symmetric neural networks, see  \cite{gerlauphawar21a}. However,  the outputs 
 obtained by these deep learning schemes only provide an approximation of the solution for a given initial distribution of the state process. 
 Hence, for another distribution $\mu_0$ of the initial state, these algorithms have to be run again. 

 In this paper, we aim to compute the minimal cost function  $v(\mu_0)$ for any $\mu_0$ $\in$ $\Pc_2(\R^d)$, and to find the optimal control, which can be searched w.l.o.g.  in the class of feedback controls, i.e., of the form 
  $\alpha_t$ $=$ $\mfa(t,X_t,\P_{X_t})$, $0\leq  t\leq T$, for some measurable function $\mfa$ on $[0,T]\times\R^d\times\Pc_2(\R^d)$. In other words, our goal is to learn 
 the value function and the optimal feedback control on the Wasserstein space.  We shall rely on a new class of neural networks,  introduced in our companion paper \cite{phawar22a}, called mean-field neural networks with input a probability measure in order to 
 approximate mappings on the Wasserstein space. 
 We then develop several numerical schemes based either on dynamic programming (DP) or stochastic maximum principle (SMP). We first propose, 
 in the spirit of \cite{gobmun05}, \cite{han2016dlapproxforscp}  a global learning of the feedback control approximated by a mean-field neural network. In the DP approach, we then propose two algorithms inspired by \cite{huretal21}: the first one learns the control by policy iteration while the second one learns sequentially the control and the value function by value iteration. In the SMP approach, we exploit the backward SDE characterization of the solution, and propose five different algorithms in line with recent methods developed in the context of standard BSDE (see 
 \cite{weinan2017deep}, \cite{hure2020deep}, \cite{germain2020deep})  that we extend to MKV BSDE  with various choices of global or local  loss functions  to be minimized in the  training of  mean-field neural networks.  
 We then provide extensive numerical experiments on three examples: a mean-field systemic risk model, a min/max linear quadratic model, and the classical mean-variance problem. We compare and discuss the advantages and drawbacks of all our algorithms.

The rest of the paper is organized as follows. We recall in Section \ref{secpreli} some  key results about the characterization of MKV control problems by DP or SMP approach, and  introduce the class of mean-field neural networks. Section \ref{secDP} presents three algorithms based on DP, while Section \ref{secBSDE} develops five algorithms based on the BSDE representation of the solution to MKV.  The performances  of all our algorithms are illustrated via three examples  in Section \ref{secnum}. Finally, we give  in 
Section \ref{seccon} some concluding remarks about the pros and cons  of the different schemes.

 \section{Preliminaries} \label{secpreli}

 \subsection{Characterization of McKean-Vlasov control}

Solution to the MKV control problem \eqref{MKVcontrol} can be characterized by dynamic programming (DP)  or maximum principle methods (see \cite{cardel19} for a detailed treatment of this topic). We recall the main results  that will be used  for designing our algorithms. 
In the DP approach, one considers  the dynamic version of problem 
\eqref{MKVcontrol} by defining the decoupled value function $V$ defined on $[0,T]\times\R^d\times\Pc_2(\R^d)$, which satisfies the backward recursion:
\begin{align}
V(t,X_t,\P_{X_t}) &= \;  \inf_{\alpha\in\Ac} \E \Big[ \int_t^{t+h} f(X_s,\P_{X_s},\alpha_s) \d s + V(t+h,X_{t+h},\P_{X_{t+h}}) \big| \Fc_t \Big], 
\end{align}
for any $t$ $\in$ $[0,T)$, $h$ $\in$ $(0,T-t]$, and starting from the terminal condition $V(T,x,\mu)$ $=$ $g(x,\mu)$, for $(x,\mu)$ $\in$ $[0,T]\times\Pc_2(\R^d)$, so that $v(\mu_0)$ $=$ $\E[V(0,X_0,\mu_0)]$.  
By sending $h$ to zero, we derive the master Bellman equation for the value function (see section 6.5.2 in \cite{cardel19})
\begin{align}
\partial_t V(t,x,\mu) + b\big(x,\mu,\hat\mra(x,\mu,\Uc(t,x,\mu),\partial_x \Uc(t,x,\mu)) \big)\cdot \partial_x V(t,x,\mu) & \\
\quad \quad + \;  \frac{1}{2} \sigma\sigma\trans(x,\mu,\hat\mra(x,\mu,\Uc(t,x,\mu),\partial_x \Uc(t,x,\mu)) \cdot \partial_{xx}^2 V(t,x,\mu) & \\
+ \;  \E_{\xi\sim\mu} \Big[ b\big(\xi,\mu,\hat\mra(\xi,\mu,\Uc(t,\xi,\mu),\partial_x \Uc(t,\xi,\mu)) \big)\cdot \partial_\mu V(t,x,\mu)(\xi)  & \\
\quad + \;  \frac{1}{2} \sigma\sigma\trans(\xi,\mu,\hat\mra(\xi,\mu,\Uc(t,\xi,\mu),\partial_x \Uc(t,\xi,\mu)) \cdot \partial_{x'}\partial_\mu V(t,x,\mu)(\xi)  \Big] & \\
\quad + \;   f\big(x,\mu,\hat\mra(x,\mu,\Uc(t,x,\mu),\partial_x \Uc(t,x,\mu)) \big)   & = \; 0,
\end{align}
for $(t,x,\mu)$ $\in$ $[0,T)\times\R^d\times\Pc_2(\R^d)$. Here $\cdot$ is the inner product in Euclidian spaces, $\trans$ is the transpose operator for a matrix, 
$x'$ $\in$ $\R^d$ $\mapsto$ $\partial_\mu V(t,x,\mu)(x')$ $\in$ $\R^d$ is the Lions derivative on $\Pc_2(\R^d)$ (see \cite{cardel19}),  the notation $\E_{\xi\sim\mu}[.]$ means that the expectation is taken w.r.t. the random variable $\xi$ distributed according to the law $\mu$, 
\begin{align} \label{defUc} 
\Uc(t,x,\mu) & = \; 
\partial_x V(t,x,\mu) + \E_{\xi\sim\mu}\big[\partial_\mu V(t,\xi,\mu)(x) \big] \\
&= \; \partial_\mu v(t,\mu)(x), \quad \mbox{ with } v(t,\mu) \; := \;  \E_{\xi\sim\mu}[V(t,\xi,\mu)],
\end{align}
and it is  assumed that  for any $(x,\mu,p,M)$ $\in$ $\R^d\times\Pc_2(\R^d)\times\R^d\times\R^{d\times d}$, there exists a minimizer 
\begin{align}
\hat\mra(x,\mu,p,M) & \in \;  \argmin_{a\in A} H(x,\mu,p,M,a),  \\
\mbox{ with } \quad H(x,\mu,p,M,a) & :=\;  b(x,\mu,a)\cdot p + \frac{1}{2} \sigma\sigma\trans(x,\mu,a)\cdot M + f(x,\mu,a), 
\end{align} 
which is Lipschitz in all its variables, so that we get an optimal feedback control given by 
\begin{align} \label{defhata} 
\mfa^\star(t,x,\mu) & = \; \hat\mra(x,\mu,\Uc(t,x,\mu),\partial_x \Uc(t,x,\mu)), \quad (t,x,\mu) \in [0,T]\times\R^d\times\Pc_2(\R^d).
\end{align}

In the case where the diffusion coefficient $\sigma(x,\mu)$ does not depend on the control variable $a$, and so $\hat\mra(x,\mu,p)$ does not depend on the variable $M$,  we have a probabilistic characterization of the solution in terms of forward-backward SDE of MKV type:  by setting 
\begin{align}
Y_t \; = \; V(t,X_t,\P_{X_t}), \quad Z_t \; = \; \sigma(X_t,\P_{X_t})\trans\partial_x V(t,X_t,\P_{X_t}), \quad 0 \leq t \leq T, 
\end{align}
it follows from It\^o's formula and Master Bellman equation that $(X,Y,Z)$ satisfies the forward-backward SDE
\begin{equation} \label{BSDEY}
\left\{
\begin{array}{ccl}
\d X_t &= & b(X_t,\P_{X_t},\hat\mra(X_t,\P_{X_t},P_t) \big) \d t + \sigma(X_t,\P_{X_t}) \d W_t, \quad 0\leq t\leq T,  \; X_0 \sim \mu_0 \\
\d Y_t & = & -  f\big(X_t,\P_{X_t},\hat\mra(X_t,\P_{X_t},P_t) \big) \d t + Z_t \cdot \d W_t, \quad 0 \leq t \leq T, \; Y_T = g(X_T,\P_{X_T}), 
\end{array}
\right.
\end{equation} 
where the pair $(P_t,M_t)_t$ $=$ $(\Uc(t,X_t,\P_{X_t}),\partial_x \Uc(t,X_t,\P_{X_t})\sigma(X_t,\P_{X_t})  )_t$ of adjoint processes, valued in $\R^d\times\R^{d\times d}$, satisfies from the Pontryagin maximum principle the backward SDE: 
\begin{equation} \label{BSDEP}
\left\{
\begin{array}{ccl}
\d P_t &= &  - \partial_x H\big(X_t,\P_{X_t},P_t,M_t,\hat\mra(X_t,\P_{X_t},P_t) \big) \d t \\
& & \quad - \; \tilde \E\Big[ \partial_\mu H\big(\tilde X_t,\P_{X_t},\tilde P_t,\tilde M_t,\hat\mra(\tilde X_t,\P_{X_t},\tilde P_t) \big)(X_t) \Big] \d t + M_t \d W_t, \quad 0 \leq t\leq T, \\
P_T &= &  \partial_x g(X_T,\P_{X_T}) + \tilde\E\big[ \partial_\mu g(\tilde X_T,\P_{X_T})(X_T) \big], 
\end{array}
\right. 
\end{equation}
where $(\tilde X,\tilde P,\tilde M)$ are independent copies of $(X,P,M)$ on $(\tilde\Omega,\tilde\Fc,\tilde\P)$. Under the assumption that $(x,\mu)$ $\in$ $\R^d\times\Pc_2(\R^d)$  $\mapsto$ $g(x,\mu)$ is convex,  $(x,\mu,a)$ $\in$ $\R^d\times\Pc_2(\R^d)\times A$ (with $A$ convex set) 
$\mapsto$  $H(x,\mu,p,M,a)$ is convex for any $(p,M)$, together with additional regularity conditions on the coefficients $b,\sigma,f,g$, it is known from \cite{carmona2015forward} that the solution to the adjoint BSDE \eqref{BSDEP}  yields an optimal control given by 
\begin{align}
\alpha_t^* &= \;  \mfa^\star(t,X_t,\P_{X_t}) \; = \; \hat\mra(X_t,\P_{X_t},P_t),  \quad 0\leq t\leq T. 
\end{align}

 We are then led to consider the generic form of MKV forward-backward $(X,\Yc,\Zc)$: 
 \begin{equation} \label{BSDEgen}
 \left\{
 \begin{array}{ccl}
 \d X_t &=& B(X_t,\P_{X_t},\Yc_t) \d t +  \sigma(X_t,\P_{X_t}) \d W_t, \quad 0\leq t\leq T,  \; X_0 \sim \mu_0,  \\
 \d\Yc_t &=& \tilde\E\big[ \Hc(X_t,\P_{X_t},\Yc_t,\Zc_t,\tilde X_t,\tilde\Yc_t,\tilde\Zc_t) \big] \d t + \Zc_t \d W_t, \quad 0 \leq t\leq T, \; \Yc_T = G(X_T,\P_{X_T}). 
 \end{array}
 \right. 
 \end{equation}

\subsection{Mean-field neural networks} 
 \label{subsec:MFNN}

\begin{figure}[H]
    \centering
    \includegraphics[scale=0.25]{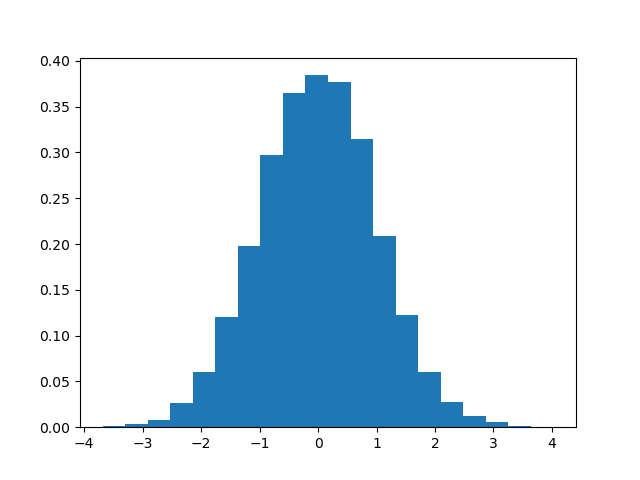}
    \caption{\small{Bin approximation of a Gaussian distribution.}}
    \label{fig:binDist}
\end{figure}

The solution to MKV control problem, i.e.,  the value function and optimal feedback control, are ma\-ppings of the state process and its probability distribution.  In order to approximate such mappings, we shall rely on mean-field neural networks introduced in our companion paper 
\cite{phawar22a}. Those are mappings 
\beqs
\Nc: \; \mu \in \Pc_2(\R^d) & \mapsto & \Nc(\mu)(\cdot) : \R^d \rightarrow \R^p,  \;\; \mbox{ with quadratic growth condition}, 
\enqs 
in one of the following forms:
\begin{itemize}
\item[(i)] {\it Bin density}:  $\Nc(\mu)(x)$ $=$ $\Phi(x,\bop^\mu)$, for $x$ $\in$ $\R^d$, $\mu$ $\in$ $\Dc_2(\R^d)$ the subset of probability measures $\mu$ in $\Pc_2(\R^d)$ which admit density functions $\mrp^\mu$  with respect to the Lebesgue measure $\lambda_d$ on $\R^d$.
Here, $\Phi$ is a standard feedforward neural network from $\R^d\times\R^K$ into $\R^p$, and $\bop^\mu$ $=$ $(p_k^\mu)_{k\in\llbracket 1,K\rrbracket}$  is the bin weight of the discrete density approximation of $\mrp^\mu$ on a fixed  bounded rectangular domain $\Kc$ of $\R^d$ divided into 
$K$ bins: $\cup_{k=1}^K{\rm Bin}(k)$ $=$ $\Kc$, of center $x_k$,   with same area size $h$ $=$ $\lambda_d(\Kc)/K$, hence given by (see Figure \ref{fig:binDist} in the case of  one dimensional  Gaussian distribution for $\mu$): 
\beqs
p^\mu_k & =&  \frac{\mrp^\mu(x_k)}{\sum_{k=1}^K \mrp^\mu(x_k) h}, \quad  k = 1,\ldots,K.  
\enqs
\item[(ii)] {\it Cylindrical}: $\Nc(\mu)(\cdot)$ $=$ $\Psi(\cdot,<\varphi,\mu>)$, where $\Psi$ is a feedforward network function (outer neural network) from $\R^d\times\R^q$ into $\R^p$, and $\varphi$ is another feedforward network function (inner neural network)  from $\R^d$ into $\R^q$ (called latent space).  
Here we denote $<\varphi,\mu>$ $:=$ $\int \varphi(x) \mu(\d x)$.  
\end{itemize}
The  relevance of mean-field neural networks is theoretically justified in \cite{phawar22a} by universal approximation theorems, and it has  been also shown how they can be trained accurately from  samples of probability measures $\mu$ $=$ $\Lc_D(\bop)$  with discrete density of bin weight $\bop$ $=$ $(p^k)_{k\in\llbracket 1,K\rrbracket}$ drawn 
randomly on $\Dc_K$ $=$ $\{\bop=(p_k)_{k\in\llbracket 1,K\rrbracket} \in \R_+^K: \sum_{k=1}^K p_k h=1\}$, and simulations of random variables $X$ $\sim$ $\mu$ by inverse transform.  Notice that for  $\mu$ $=$ $\Lc_D(\bop)$, we have $\bop^\mu$ $=$ $\bop$, and so the bin density network at such $\mu$ is equal to  
$\Nc(\mu)(.)$ $=$ $\Phi(.,\bop)$. On the other hand, for any cylindrical function $F$ of the measure  of the form $F(\mu)$ $=$ $\Psi(<\varphi,\mu>)$,  we can compute it approximately from samples $X^{(n)}$, $n$ $=$ $1,\ldots,N$, of $\mu$ by:
$F(\mu)$ $\simeq$ $\Psi\big( \frac{1}{N} \sum_{n=1}^N \varphi(X^{(n)}) \big)$. This is the case  in particular for cylindrical neural network. 

As described in \cite{phawar22a}, exploring the space of probability measures is crucial for both neural networks. In both cases, we employ the bins method to generate samples of probability measures for training mean-field neural networks. The algorithm used to generate these measures is outlined in \cite{phawar22a} and is currently limited to dimension one. Consequently, all numerical tests conducted in the article are confined to dimension one. However, it is possible to handle cases in dimension two by employing a different algorithm proposed in \cite{warin2023quantile}. 
In all subsequent algorithms, the proper selection of the domain $\Kc$ is crucial, particularly for the bins method. When the support of the distribution is unknown, an iterative procedure becomes necessary.
Two algorithms can be implemented as follows:
\begin{itemize}
\item[1.] First algorithm: (i) Initially, make an initial guess of the support.
(ii) Once the resolution is obtained, verify that the generated distribution's support is primarily contained within $\Kc$, sufficiently far from its boundary. (iii) If the support is not mainly within $\Kc$, adapt the size and center of $\Kc$ accordingly.
\item[2.] Second algorithm: (i) Use a very large $\Kc$ during the first iteration to locate the domain of importance, employing a coarse resolution. (ii) In the subsequent calculation, reduce the size of $\Kc$ to achieve an accurate resolution.
\end{itemize}

\section{Dynamic programming-based algorithms} \label{secDP} 

We consider a time discretization of the MKV control problem by fixing a time grid  $\Tc$ $=$ $\{t_i=i\Delta t: i=0,\ldots,N_T\}$, 
with $\Delta t$ $=$ $T/N_T$, and introducing the corresponding 
mean-field Markov decision process: minimize over feedback controls $\mfa$ on $\Tc\times\R^d\times\Pc_2(\R^d)$ the cost functional
\beqs
J_{N_T}(\mfa) &=& \E \Big[ \sum_{i=0}^{N_T-1} f(X_i,\mu_{i},\mfa(t_i,X_i,\mu_{i})) \Delta t + g(X_{N_T},\mu_{{N_T}}) \Big], 
\enqs
where 
\begin{align}
X_{i+1} &= \;  X_i + b(X_i,\mu_{i},\mfa(t_i,X_i,\mu_{i})) \Delta t + \sigma(X_i,\mu_{i},\mfa(t_i,X_i,\mu_{i})) \Delta W_i,  \\
& =: \;  F_{\Delta t}(X_i,\mu_{i},\mfa(t_i,X_i,\mu_{i}),\Delta W_i),  \quad   i=0,\ldots,N_T -1, \; X_0 \sim \mu_0,
\end{align}
with  $\Delta W_i$ $:=$ $W_{t_{i+1}}-W_{t_i}$, and $\mu_i$ $=$ $\P_{X_i}$ denotes the law of $X_i$. 

We present two classes of algorithms. The first one is learning the control by a single optimization but allows us  to compute the solution of the problem \eqref{MKVcontrol}  and therefore the solution of the correspon\-ding master Bellman  equation only at time  $t=0$ for all distributions $\mu_0$.
The second class with two other algorithms solves  $N_T$ local optimization  problems, and allows us to  compute the master equation at all dates for all distributions.

\subsection{Global learning on control} \label{sec:globalcon}

In the spirit of the method introduced in \cite{gobmun05}, \cite{han2016dlapproxforscp},  which does not actually rely on dynamic programming, 
we replace feedback controls  by  time-dependent mean-field neural networks $\Nc(t,\mu)(x)$ valued in $A$ $\subset$ $\R^m$, with input $t$ $\in$ $[0,T]$, $\mu$ $\in$ $\Pc_2(\R^d)$, and $x$ $\in$ $\R^d$, and minimize over the parameters $\theta$ of this 
mean-field neural network $\Nc$ $=$ $\Nc_\theta$ the global cost function
\beqs
\mrJ(\theta) &=& \E \Big[  \sum_{i=0}^{N_T-1} f(X_i,\mu_{i},\Nc_\theta(t_i,\mu_{i})(X_i)) \Delta t + g(X_{N_T},\mu_{{N_T}}) \Big], 
\enqs
with 
\beqs
X_{i+1} &=&  F_{\Delta t}(X_i,\mu_{i},\Nc_\theta(t_i,\mu_{i})(X_i),\Delta W_i),  \quad   i=0,\ldots,N_T -1, \; X_0 \sim \mu_0. 
\enqs
 In practice, for $i$ $=$ $1,\ldots,N_T$, $\mu_{i}$ has to be estimated/approximated from samples of $X_{i}$, and this is done as follows. We use a training batch of $M$  probability measures $\mu_0^{(m)}$ $=$ $\Lc_D(\bop^{(m)})$ in $\Dc_2(\R^d)$ from samples 
 $\bop^{(m)}$ $=$ $(p_k^{(m)})_{k\in\llbracket 1,K\rrbracket}$, $m$ $=$ $1,\ldots,M$, in $\Dc_K$. Then, for each $m$, we sample $X_0^{(m),(n)}$, $n$ $=$ $1,\ldots,N$, from $\mu_0^{(m)}$, and for $i$ $=$ $0,\ldots,N_T-1$,  $X_{i+1}^{(m),(n)}$, $n$ $=$ $1,\ldots,N$ are sampled as 
\beqs
X_{i+1}^{(m),(n)} &=& F_{\Delta t}(X_i^{(m),(n)},\hat\mu_i^{(m)},\Nc_\theta(t_i,\hat\mu^{(m)}_{i})(X_i^{(m),(n)}),\Delta W_{i}^{(m),(n)}),   
\enqs
with $\hat\mu_i^{(m)}$ $=$ $\Lc_D(\hat\bop_i^{(m)})$, $\hat\bop_0^{(m)}$ $=$ $\bop^{(m)}$, and $\hat\bop_i^{(m)}$ $=$ $(\hat p_{i,k}^{(m)})_{k \in \llbracket 1,K\rrbracket}$ are the estimated density weights in $\Dc_K$ of  $X_{i}^{(m),(n)}$, $i$ $=$ $1,\ldots,N_T$ (truncated on $\Kc$), namely: 
\beqs
\hat p_{i,k}^{(m)} &=& \frac{ \#  \{ n \in \llbracket 1,N\rrbracket:  {\rm Proj}_\Kc(X_{i}^{(m),(n)}) \in \mbox{Bin}(k) \} }{Nh}, \quad k=1,\ldots,K,
\enqs
where ${\rm Proj}_{\Kc}(.)$ is the projection on $\Kc$.  The cost function is then approximated by 
\beqs
\mrJ_{M,N}(\theta) & = & \frac{1}{MN} \sum_{m=1}^M \Big[  \sum_{n=1}^N  
\sum_{i=0}^{N_T-1}  f\big(X_i^{(m),(n)},\hat\mu_i^{(m)},\Nc_\theta(t_i,\hat\mu_{i}^{(m)})(X_i^{(m),(n)})  \big) \Delta t + g(X_{N_T}^{(m),(n)},\hat\mu_{N_T}^{(m)}) \Big].
\enqs

\vspace{1mm}

The pseudo-code using a gradient descent method  is described in Algorithm \ref{algo:globalcont}. 

\vspace{1mm}

\begin{algorithm}[H] 
{\bf Input data:} A time-dependent mean-field neural network $\Nc_\theta(t,\mu)(x)$.  

{\bf Initialization:} learning rate $\gamma$ and parameters $\theta$ 

\For{each epoch} 
{
Generate a batch  of $M$ initial distributions $\mu_0^{(m)}$, $m$ $=$ $1,\ldots,M$ \;
\For{$m=1,\ldots,M$}
{Generate Brownian increments $\Delta W_i^{(m),(n)}$, $i$ $=$ $0,\ldots,N_T-1$, $n$ $=$ $1,\ldots,N$ \;
Compute sample trajectories $X_0^{(m),(n)}$, $X_i^{(m),(n)}$, $n$ $=$ $1,\ldots,N$, and estimate $\hat\mu_i^{(m)}$, $i$ $=$ $1,\ldots,N_T$,  
}
Compute the batch cost $\mrJ_{M,N}(\theta)$ and its gradient $\nabla_\theta \mrJ_{M,N}(\theta)$ \;
Update $\theta$ $\leftarrow$ $\theta$ $-$ $\gamma \nabla_\theta \mrJ_{M,N}(\theta)$ \;
}
{\bf Return:} The set of optimized parameters $\theta^*$. 
\caption{Global learning  on the control} \label{algo:globalcont} 
\end{algorithm}

\vspace{1mm}

The global algorithms that directly minimize the objective function have demonstrated effectiveness in practice, even without having a theoretical convergence proof. 
The output of this global algorithm is an approximation of the optimal feedback control at initial time $t_0$ $=$ $0$ by a mean-field neural network  $\Nc_{\theta^*}(t_0,.)$, and yields an approximation of the optimal control at other times $t_i$, $i$ $=$ $1,\ldots,N_T-1$, 
by mean-field neural networks $\Nc_{\theta^*}(t_i,\mu_i)(X_i)$ along the law  $\mu_i$, and the  state $X_i$ explored during the learning algorithm. 
The value function can then be estimated at initial time $t_0$  by regression as follows: we approximate the initial value function by a mean-field neural neural network $\vartheta_\eta(\mu)(x)$ valued in $\R$, and minimize over the parameters $\eta$ of this neural network the quadratic loss function 
 \beqs
 \E  
 \Big|  \sum_{i=0}^{N_T-1} f(X_i,\mu_{i},\Nc_{\theta^*}(t_i,\mu_{i})(X_i)) \Delta t + g(X_{N_T},\mu_{{N_T}})  - \vartheta_\eta(\mu_0)(X_0)  \Big|^2, 
 \enqs
 where 
\beqs
X_{i+1} &=&  F_{\Delta t}(X_i,\mu_{i},\Nc_{\theta^*}(t_i,\mu_{i})(X_i),\Delta W_i),  \quad   i=0,\ldots,N_T -1, \; X_0 \sim \mu_0. 
\enqs
When using the global algorithm and the cylinder network, there is no need to estimate the support of the distribution. The parameter $\Kc$ is solely used to generate probability distributions at time $0$, and its selection is based on ensuring that the initial distribution of $X_0$ primarily concentrates its mass within $\Kc$.
On the other hand, when employing the bin method, it is necessary to monitor the generated distribution and verify that its support is predominantly contained within $\Kc$. If this is not the case, the size of $\Kc$ should be adjusted using the procedure suggested in Section \ref{subsec:MFNN}.

\subsection{Control learning by policy iteration} \label{seccontNN}

Our next algorithm is inspired by the method in \cite{huretal21}, which is a combination of the global algorithm on control and dynamic programming. We replace at any time $t_i$, $i$ $=$ $0,\ldots,N_T-1$, feedback controls by mean-field neural networks 
$\Nc_{\theta_i}$ with parameter $\theta_i$,  and proceed by backward induction for computing approximate optimal controls: for $i$ $=$ $N_T-1,\ldots,0$, keep track of the approximate optimal feedback  controls $\Nc_{\theta_j^*}$, $j$ $=$ $i+1,\ldots,N_T-1$, and minimize over 
$\theta_i$ the cost function: 
\beqs
\mrJ^i(\theta_i) &=& \E \Big[ f(X_i,\mu_i,\Nc_{\theta_i}(\mu_i)(X_i)) \Delta t +  \sum_{j=i+1}^{N_T-1}  f(X_j,\mu_j,\Nc_{\theta_j^*}(\mu_j)(X_j)) \Delta t + g(X_{N_T},\mu_{N_T}) \Big],
\enqs 
(with the convention that the above sum over $j$ is empty when $i$ $=$ $N_T-1$) where 
\begin{equation} \label{XijcontNN}
\begin{cases}
X_{i+1} \; = \;   F_{\Delta t}(X_i,\mu_{i},\Nc_{\theta_i}(\mu_{i})(X_i),\Delta W_i), \quad X_i \sim \mu_i,  \\
X_{j+1} \; = \;  F_{\Delta t}(X_j,\mu_{j},\Nc_{\theta_j^*}(\mu_{j})(X_j),\Delta W_j), \quad j=i+1,\ldots,N_T -1. 
\end{cases}
\end{equation} 
In the practical implementation, the  cost function $\mrJ^i(.)$ is approximately computed from a training of $M$ probability measures  $\mu_i^{(m)}$ $=$ $\Lc_D(\bop_i^{(m)})$ in $\Dc_2(\R^d)$ with samples $\bop_i^{(m)}$ $=$ $(p_{i,k}^{(m)})_{k\in\llbracket 1,K\rrbracket}$, $m$ $=$ $1,\ldots,M$, in $\Dc_K$. 
For each batch $m$, one then computes $N$ samples $X_i^{(m),(n)}$ $\sim$ $\mu_i^{(m)}$,  $X_j^{(m),(n)}$, $j$ $=$ $i+1,\ldots,N_T-1$, $n$ $=$ $1,\ldots,N$, according to \eqref{XijcontNN} with estimated probability measures  $\hat\mu_j^{(m)}$ $=$ $\Lc_D(\hat\bop_j^{(m)})$, as in Section \ref{sec:globalcon}, and thus approximate 
the local cost function by 
\begin{align}
\mrJ^i_{M,N}(\theta_i) &=\;    \frac{1}{MN} \sum_{m=1}^M \sum_{n=1}^N \Big[ f(X_i^{(m),(n)},\mu_i^{(m)},\Nc_{\theta_i}(\mu_i^{(m)})(X_i^{(m),(n)})) \Delta t \\
&  \quad  \quad \quad +   \sum_{j=i+1}^{N_T-1}  f(X_j^{(m),(n)},\hat\mu_j^{(m)},\Nc_{\theta_j^*}(\hat\mu_j^{(m)})(X_j^{(m),(n)})) \Delta t
+ g(X_{N_T}^{(m),(n)},\hat\mu^{(m)}_{N_T}) \Big]. 
\end{align}

\vspace{1mm}

The pseudo-code is described in Algorithm \ref{alg:fControl}. 

\vspace{1mm}

\begin{algorithm}[H]
{\bf Input data:} Mean-field neural networks $\Nc_{\theta_i}$  \;
\For{$i= N_T -1,\ldots,0$}
 { {\it Initialization:} learning rate $\gamma$ and parameters $\theta_i$   \; 
  \For{each epoch} 
  {
    Generate a batch  of $M$ distributions $\mu_i^{(m)}$, $m$ $=$ $1,\ldots,M$ \;  
    \For{$m$ $=$ $1,\ldots,M$}
    {
    Generate Brownian increments $\Delta W_k^{(m),(n)}$, $k$ $=$ $i,\ldots,N_T-1$, $n$ $=$ $1,\ldots,N$ \;
    Compute sample trajectories $X_i^{(m),(n)}$, $X_j^{(m),(n)}$, $n$ $=$ $1,\ldots,N$, and estimate $\hat\mu_j^{(m)}$, $j$ $=$ $i+1,\ldots,N_T$,  
    }
Compute the batch cost $\mrJ^i_{M,N}(\theta_i)$ and its gradient $\nabla_\theta \mrJ^i_{M,N}(\theta_i)$ \;
Update $\theta_i$ $\leftarrow$ $\theta_i$ $-$ $\gamma \nabla_\theta \mrJ^i_{M,N}(\theta_i)$ \;
}
$\theta_i^* = \theta_i$
}
{\bf Return:}  Optimized parameters $\theta_i^*$, $i=0, \ldots, N_T-1$. 
      \caption{Learning by policy iteration \label{alg:fControl}}
\end{algorithm}

\vspace{2mm}

The output of this  algorithm is an approximation of the optimal feedback control at any time $t_i$ by a mean-field neural network   $\Nc_{\theta_i^*}$, $i$ $=$ $0,\ldots,N_T-1$. 
 The value function can then be estimated at any time $t_i$  by regression as follows: we approximate the value function at time $t_i$ by a mean-field neural neural network $\vartheta_{\eta_i}(\mu)(x)$ valued in $\R$, and minimize over the parameters $\eta_i$ of this neural network the quadratic loss function 
 \begin{align}  \label{regressV} 
 \E   \Big|  \sum_{j=i}^{N_T-1} f(X_j,\mu_{j},\Nc_{\theta_j^*}(\mu_{j})(X_j)) \Delta t + g(X_{N_T},\mu_{{N_T}})  - \vartheta_{\eta_i}(\mu_i)(X_i)  \Big|^2, 
 \end{align} 
 where 
\beqs
X_{j+1} &=&  F_{\Delta t}(X_j,\mu_{j},\Nc_{\theta_j^*}(\mu_{j})(X_j),\Delta W_j),  \quad   j=i,\ldots,N_T -1, \; X_i \sim \mu_i. 
\enqs
In a backward algorithm, having a good estimate of the support of the distribution being tested is crucial at each time step $i$. This estimate helps in efficiently sampling the distribution in areas of interest. If the support is unknown, an iterative procedure, such as the one proposed in Section \ref{subsec:MFNN}, needs to be implemented to gradually refine the estimation of the support.

\subsection{Control learning by value iteration}

The two previous algorithms provide low bias estimates of the learnt controls, but in general high-variance estimate due to this cumulated sum over the cost functions. 
Moreover, these algorithms are very memory demanding  as, at each epoch, all the $N$ trajectories for the  $M$ distributions have to be generated for the $O(N_T)$ time values  and stored. 
To circumvent this possible variance issue, we propose an alternate algorithm of actor-critic type, similarly as in \cite{huretal21} (called there hybrid algorithm), where the feedback control and value function are learnt  sequentially. We are given a family of mean-field neural networks $\Nc_{\theta_i}$ and $\vartheta_{\eta_i}$, $i$ $=$ $0,\ldots,N_T-1$, for the approximation of the feedback control (actor) and value function (critic).  
We proceed by backward induction as follows: starting from $\vartheta_{N_T}^*(\mu)(x)$ $=$ $g(x,\mu)$,  we minimize over $\theta_i$, for $i$ $=$ $N_T-1,\ldots,0$, the  cost function
\beqs
\mrJ^i(\theta_i) &=& \E \Big[ f(X_i,\mu_i,\Nc_{\theta_i}(\mu_i)(X_i)) \Delta t + \vartheta_{i+1}^*(\mu_{i+1})(X_{i+1}) \Big],
\enqs
where 
\begin{align} \label{Xi+1} 
X_{i+1} &= \;   F_{\Delta t}(X_i,\mu_{i},\Nc_\theta(t_i,\mu_{i})(X_i),\Delta W_i), \; X_i \sim \mu_i, 
\end{align} 
update $\theta_i^*$ as the resulting optimal parameter, then minimize over $\eta_i$ the quadratic loss function
\beqs
L^i(\eta_i) &=& \E \Big| f(X_i,\mu_i,\Nc_{\theta^*_i}(\mu_i)(X_i)) \Delta t + \vartheta_{i+1}^*(\mu_{i+1})(X_{i+1})  -  \vartheta_{\eta_i}(\mu_i)(X_i) \Big|^2, 
\enqs
update $\eta_i^*$ as the resulting optimal parameter, and set $\vartheta_i^*$ $=$  $\vartheta_{\eta_i^*}$.   
Again, in the practical implementation, we use a training of $M$ probability measures  $\mu_i^{(m)}$ $=$ $\Lc_D(\bop_i^{(m)})$ in $\Dc_2(\R^d)$ with samples $\bop_i^{(m)}$ $=$ $(p_{i,k}^{(m)})_{k\in\llbracket 1,K\rrbracket}$, $m$ $=$ $1,\ldots,M$, in $\Dc_K$. 
For each batch $m$, one then computes $N$ samples $X_i^{(m),(n)}$ $\sim$ $\mu_i^{(m)}$,  $X_{i+1}^{(m),(n)}$ according to \eqref{Xi+1} with estimated probability measure  $\hat\mu_{i+1}^{(m)}$ $=$ $\Lc_D(\hat\bop_{i+1}^{(m)})$, as in Section \ref{sec:globalcon}, and approximate the  function $\mrJ^i$ by  
\begin{align}
\mrJ^i_{M,N}(\theta_i) &=\;    \frac{1}{MN} \sum_{m=1}^M \sum_{n=1}^N \Big[ f(X_i^{(m),(n)},\mu_i^{(m)},\Nc_{\theta_i}(\mu_i^{(m)})(X_i^{(m),(n)})) \Delta t +  \vartheta_{i+1}^*(\mu_{i+1}^{(m)})(X_{i+1}^{(m),(n)}) \Big],
\end{align}
while similarly  the second loss function $L^i$  is approximated by
\begin{align}
L^i_{M,N}(\eta_i) =&  \frac{1}{MN} \sum_{m=1}^M \sum_{n=1}^N  \Big| f(X_i^{(m),(n)},\mu_i^{(m)},\Nc_{\theta^*_i}(\mu_i^{(m)})(X_i^{(m),(n)})) \Delta t + \vartheta_{i+1}^*(\mu_{i+1}^{(m)})(X_{i+1}^{(m),(n)})  \\
 & \hspace{2.5cm} - \;   \vartheta_{\eta_i}(\mu_i^{(m)})(X_i^{(m),(n)}) \Big|^2. 
\end{align}

\vspace{2mm}

The pseudo-code is described in Algorithm \ref{alg:LocalfVal}. 

\vspace{1mm}

\begin{algorithm}[H]
{\bf Input data:} Mean-field neural networks $\Nc_{\theta_i}$, $\vartheta_{\eta_i}$,  $i$ $=$ $0,\ldots,N_T-1$ \;  

{\bf Initialization}: $\vartheta^*_{{N_T}}(\mu)(x)$ $=$ $g(x,\mu)$;

\For{$i= N_T -1,\ldots,0$}
{{\it Initialization:} learning rates $\gamma_C$, $\gamma_V$ and parameters $\theta_i$, $\eta_i$ \;

 \For{each epoch} 
{
Generate a batch  of $M$ distributions $\mu_i^{(m)}$, $m$ $=$ $1,\ldots,M$ \;  
 \For{each batch $m$}
 {
    Generate Brownian increments $\Delta W_i^{(m),(n)}$,  $n$ $=$ $1,\ldots,N$ \;
Compute samples  $X_i^{(m),(n)}$, $X_{i+1}^{(m),(n)}$, $n$ $=$ $1,\ldots,N$, and estimate $\hat\mu_{i+1}^{(m)}$ 
    }
Compute the batch cost $\mrJ^i_{M,N}(\theta_i)$ and its gradient $\nabla_\theta \mrJ^i_{M,N}(\theta_i)$ \;
Update $\theta_i$ $\leftarrow$ $\theta_i$ $-$ $\gamma \nabla_\theta \mrJ^i_{M,N}(\theta_i)$ \;
}
Store optimized parameter $\theta_i^*$ \;
\For{each epoch} 
{
Generate a batch  of $M$ distributions $\mu_i^{(m)}$, $m$ $=$ $1,\ldots,M$ \;  
 \For{each batch $m$}
    {
    Generate Brownian increments $\Delta W_i^{(m),(n)}$,  $n$ $=$ $1,\ldots,N$ \;
    Compute samples  $X_i^{(m),(n)}$, $X_{i+1}^{(m),(n)}$, $n$ $=$ $1,\ldots,N$, and estimate $\hat\mu_{i+1}^{(m)}$ 
    }
Compute the batch cost $L^i_{M,N}(\eta_i)$ and its gradient $\nabla_\eta L^i_{M,N}(\eta_i)$ \;
Update $\eta_i$ $\leftarrow$ $\eta_i$ $-$ $\gamma \nabla_\eta L^i_{M,N}(\eta_i)$ \;
}
$\vartheta^*_i = \vartheta_{\eta_i^*}$ 
}
{\bf Return:}  The optimized parameters $\theta_i^*$, $\eta_i^*$, $i=0, \ldots, N_T-1$
    \caption{Actor/critic algorithm: learning by value iteration  \label{alg:LocalfVal}}
\end{algorithm}

\vspace{1mm}

The output of this  algorithm is an approximation of the optimal feedback control and value function at any time $t_i$ by  mean-field neural networks   $\Nc_{\theta_i^*}$, and $\vartheta_{\eta_i^*}$,  $i$ $=$ $0,\ldots,N_T-1$. \\
Since the resolution is performed in a backward manner, when the support of the distribution is unknown, it becomes necessary to employ an iterative algorithm, as described in Section \ref{subsec:MFNN}, to explore and estimate the distributions of interest.

\section{Backward SDE-based algorithms} \label{secBSDE}

We start from  the time discretization of the MKV forward-backward SDE \eqref{BSDEgen} that characterizes the solution to the MKV control problem: 
\begin{equation}
\begin{cases}
X_{i+1}  =  X_i  + B(X_i,\mu_i,\Yc_i) \Delta t + \sigma(X_i,\mu_i) \Delta W_i, \quad i=0,\ldots,N_T-1, \; X_0 \sim \mu_0, \\
\Yc_{i+1} = \Yc_i + \tilde\E\big[\Hc(X_i,\mu_i,\Yc_i,\Zc_i,\tilde X_i,\tilde\Yc_i,\tilde\Zc_i)\big] \Delta t + \Zc_i \Delta W_i, \quad i=0,\ldots,N_T-1, \; \Yc_{N_T} = G(X_{N_T},\mu_{N_T}). 
\end{cases}
\end{equation}

This system of equations corresponds to the resolution of the system of equations \eqref{BSDEY}, \eqref{BSDEP}. Note that in fact, $(X_t,P_t)$ is independent of $Y_t$. Then the resolution is achieved by calculating the optimal control solving $(X_t,P_t)$ for $t \le T$. The estimation of $Y_t$ is  achieved by using the optimal control with a simple forward simulation and by taking the expectation of  $Y_t$  in equation \eqref{BSDEY}:
\begin{align*}
   Y_t &= \;   \E\Big [ \int_t^T ] f\big(X_s,\P_{X_s},\hat\mra(X_s,\P_{X_s},P_s) \big) \d s +  g(X_T,\P_{X_T}) \big| \Fc_t  \Big].   
\end{align*}

\subsection{Local algorithms}


We adapt the deep backward scheme in \cite{hure2020deep} to our context. 
We are given a family of mean-field neural networks $\Yc_{\theta_i}(\mu)(x),\Zc_{\theta_i}(\mu)(x)$, $i$ $=$ $0,\ldots,N_T-1$ (by misuse of notation, we also denote by $\Yc$ and $\Zc$  the neural networks for the approximation of the pair component of the MKV BSDE), and proceed by backward induction as follows: 
starting from $\Yc_{N_T}^*(\mu)(x)$ $=$ $G(x,\mu)$,  we minimize over $\theta_i$, for $i$ $=$ $N_T-1,\ldots,0$, 
the loss function 
\beqs
L^i(\theta_i) &=& \E \Big| \Yc^*_{i+1}(\mu_{i+1})(X_{i+1}) -  \Yc_{\theta_i}(\mu_i)(X_i) -  \Zc_{\theta_i}(\mu_i)(X_i) \Delta W_i \\
& & \quad \quad - \;   \tilde\E\big[\Hc(X_i,\mu_i,\Yc_{\theta_i}(\mu_i)(X_i),\Zc_{\theta_i}(\mu_i)(X_i),\tilde X_i,\Yc_{\theta_i}(\mu_i)(\tilde X_i),\Zc_{\theta_i}(\mu_i)(\tilde X_i))\big] \Delta t  \Big|^2,
\enqs
where 
\begin{align} \label{XBSDEtheta}
X_{i+1} &= \;  X_i  + B(X_i,\mu_i,\Yc_{\theta_i}(\mu_i)(X_i)) \Delta t + \sigma(X_i,\mu_i) \Delta W_i, \quad X_i \sim \mu_i, 
\end{align} 
update $\theta_i^*$ as the resulting optimal parameter, and set $\Yc_i^*$ $=$ $\Yc_{\theta_i^*}$.  In the practical implementation, we use a training of $M$ probability measures  $\mu_i^{(m)}$ $=$ $\Lc_D(\bop_i^{(m)})$ in $\Dc_2(\R^d)$ with samples $\bop_i^{(m)}$ $=$ $(p_{i,k}^{(m)})_{k\in\llbracket 1,K\rrbracket}$, $m$ $=$ $1,\ldots,M$, in $\Dc_K$. 
For each batch $m$, one then computes $N$ independent samples $X_i^{(m),(n)}$,  $\tilde X_i^{(m),(n)}$ $\sim$ $\mu_i^{(m)}$,  $n$ $=$ $1,\ldots,N$, $X_{i+1}^{(m),(n)}$ according to \eqref{XBSDEtheta} with estimated probability measure  $\hat\mu_{i+1}^{(m)}$  as in Section \ref{sec:globalcon}, and approximate the loss function by 
\beqs
L_{M,N}^i(\theta_i) &=& \frac{1}{MN} \sum_{m=1}^M\sum_{n=1}^N  \Big| \Yc^*_{i+1}(\hat\mu_{i+1}^{(m)})(X_{i+1}^{(m),(n)}) -  \Yc_{\theta_i}(\mu_i^{(m)})(X_i^{(m),(n)}) -  \Zc_{\theta_i}(\mu_i^{(m)})(X_i^{(m),(n)}) \Delta W_i  \\
& &  \quad \quad  - \; \frac{\Delta t}{N} \sum_{n'=1}^N \Hc(X_i^{(m),(n)},\mu_i^{(m)},\Yc_{\theta_i}(\mu_i^{(m)})(X_i^{(m),(n)}),\Zc_{\theta_i}(\mu_i^{(m)})(X_i^{(m),(n)}),\\
& & \hspace{3cm}  \tilde X_i^{(m),(n')},\Yc_{\theta_i}(\mu_i^{(m)})(\tilde X_i^{(m),(n')}),\Zc_{\theta_i}(\mu_i^{(m)})(\tilde X_i^{(m),(n')}))   \Big|^2. 
\enqs

\vspace{2mm}

The pseudo-code is described in Algorithm \ref{alg:BSDElocal1}. It  is in the spirit of the actor/critic algorithm \ref{alg:LocalfVal}, but now  $\Yc$ and $\Zc$ are learnt simultaneously.

\vspace{1mm}

\begin{algorithm}[H]
{\bf Input data:} Mean-field neural networks $\Yc_{\theta_i}$, $\Zc_{\theta_i}$ \;  

{\bf Initialization}: $\Yc^*_{{N_T}}(\mu)(x)$ $=$ $G(x,\mu)$;

\For{$i= N_T -1,\ldots,0$}
{{\it Initialization:} learning rate $\gamma$,  and parameter $\theta_i$ \;
\For{each epoch} 
{
Generate a batch  of $M$ distributions $\mu_i^{(m)}$, $m$ $=$ $1,\ldots,M$ \;  
 \For{each batch $m$}
    {
    Generate Brownian increments  $\Delta W_i^{(m),(n)}$, $n$ $=$ $1,\ldots,N$ \;
    Compute samples  $X_i^{(m),(n)}$, $\tilde X_i^{(m),(n)}$, $X_{i+1}^{(m),(n)}$, $n$ $=$ $1,\ldots,N$, and estimate $\hat\mu_{i+1}^{(m)}$ 
    }
Compute the batch loss $L^i_{M,N}(\theta_i)$ and its gradient $\nabla_\theta L^i_{M,N}(\theta_i)$ \;
Update $\theta_i$ $\leftarrow$ $\theta_i$ $-$ $\gamma \nabla_\theta L^i_{M,N}(\theta_i)$; 
}
$\Yc_i^*$ $=$ $\Yc_{\theta_i^*}$
}
{\bf Return:}  The set of optimized parameters $\theta_i^*$, $i=0, \ldots, N_T-1$
    \caption{Deep backward algorithm   \label{alg:BSDElocal1}}
\end{algorithm}

 \vspace{5mm}

 We also propose a multi-step version of the above algorithm following the idea in \cite{germain2020deep}, and in the spirit of the policy iteration in Section \ref{seccontNN}. 
 We proceed by backward induction for $i$ $=$ $N_T-1,\ldots,0$, by keeping track of the approximate optimal mean-field neural networks $\Yc_j^*$, $\Zc_j^*$, $j$ $=$ $i+1,\ldots,N_T-1$, and minimize over 
 $\theta_i$ the loss function
 \beqs
\tilde  L^i(\theta_i) &=& \E\Big| G(X_{N_T},\mu_{N_T}) -   \sum_{j=i+1}^{N_T-1}  \Zc^*_{j}(\mu_j)(X_j) \Delta W_j -  \Zc_{\theta_i}(\mu_i)(X_i) \Delta W_i  -  \Yc_{\theta_i}(\mu_i)(X_i)    \\
 & & \quad -  \; \sum_{j=i+1}^{N_T-1}  \tilde\E\big[\Hc(X_j,\mu_j,\Yc^*_{j}(\mu_j)(X_j),\Zc^*_{j}(\mu_j)(X_j),\tilde X_j,\Yc^*_{j}(\mu_j)(\tilde X_j),\Zc^*_{j}(\mu_j)(\tilde X_j))\big] \Delta t \\
& &  \quad - \;   \tilde\E\big[\Hc(X_i,\mu_i,\Yc_{\theta_i}(\mu_i)(X_i),\Zc_{\theta_i}(\mu_i)(X_i),\tilde X_i,\Yc_{\theta_i}(\mu_i)(\tilde X_i),\Zc_{\theta_i}(\mu_i)(\tilde X_i))\big] \Delta t  \Big|^2,
 \enqs
where    
\begin{equation} \label{Xijmultistep}
\begin{cases}
X_{i+1} \; = \;  X_i  + B(X_i,\mu_i,\Yc_{\theta_i}(\mu_i)(X_i)) \Delta t + \sigma(X_i,\mu_i) \Delta W_i, \quad X_i \sim \mu_i, \\
X_{j+1} \; = \;  X_j  + B(X_j,\mu_j,\Yc^*_{j}(\mu_j)(X_j)) \Delta t + \sigma(X_j,\mu_j) \Delta W_j,  \quad j=i+1,\ldots,N_T -1. 
\end{cases}
\end{equation}    
In the practical implementation, we use a training of $M$ probability measures  $\mu_i^{(m)}$, $m$ $=$ $1,\ldots,M$, and for each batch $m$, one then computes $N$ samples $X_i^{(m),(n)}$,$\tilde X_i^{(m),(n)}$ $\sim$ $\mu_i^{(m)}$,  
$X_j^{(m),(n)}$, $\tilde X_j^{(m),(n)}$, $j$ $=$ $i+1,\ldots,N_T-1$, according to \eqref{Xijmultistep} with estimated probability measures  $\hat\mu_j^{(m)}$ $=$ $\Lc_D(\hat\bop_j^{(m)})$, as in Section \ref{sec:globalcon},   and approximate the loss function by $\tilde L^i_{M,N}(\theta_i)$, $i$ $=$ $0,\ldots,N_T-1$.
\vspace{2mm}

The pseudo-code is described in Algorithm \ref{alg:BSDElocal2}.  

\vspace{1mm}

\begin{algorithm}[H]
{\bf Input data:} Mean-field neural networks $\Yc_{\theta_i}$, $\Zc_{\theta_i}$, and  Brownian increments $\Delta W_i$, $i$ $=$ $0,\ldots,N_T-1$ \;  


\For{$i= N_T -1,\ldots,0$}
{{\it Initialization:} learning rate $\gamma$,  and parameter $\theta_i$ \;
\For{each epoch} 
{
Generate a batch  of $M$ distributions $\mu_i^{(m)}$, $m$ $=$ $1,\ldots,M$ \;  
 \For{each batch $m$}
    {
    Generate Brownian increments $\Delta W_k^{(m),(n)}$, $ \tilde \Delta W_k^{(m),(n)}$, $k$ $=$ $i,\ldots,N_T-1$, $n$ $=$ $1,\ldots,N$ \;
    Compute samples  $X_i^{(m),(n)}$, $\tilde X_i^{(m),(n)}$, $X_{j}^{(m),(n)}$, $\tilde X_j^{(m),(n)}$, $n$ $=$ $1,\ldots,N$, and estimate $\hat\mu_{j}^{(m)}$, $j$ $=$ $i+1,\ldots,N_T$ 
    }
 Compute the batch loss $\tilde L^i_{M,N}(\theta_i)$ and its gradient $\nabla_\theta \tilde L^i_{M,N}(\theta_i)$ \;
Update $\theta_i$ $\leftarrow$ $\theta_i$ $-$ $\gamma \nabla_\theta \tilde L^i_{M,N}(\theta_i)$; 
}
$\Yc_i^*$ $=$ $\Yc_{\theta_i^*}$, $\Zc_i^*$ $=$ $\Zc_{\theta_i^*}$
}
{\bf Return:}   The set $\Yc_i^*$ $=$ $\Yc_{\theta_i^*}$, $\Zc_i^*$ $=$ $\Zc_{\theta_i^*}$, $i=0, \ldots,N_T-1$
    \caption{Deep backward multi-step algorithm   \label{alg:BSDElocal2}}
\end{algorithm}
 
 \vspace{1mm}

The output of these two algorithms  \ref{alg:BSDElocal1} and \ref{alg:BSDElocal2} yields in particular  an approximation of the function $\Uc$ in \eqref{defUc} by the mean-field neural network $\Yc_i^*$ at any time $t_i$, hence an approximation of the optimal feedback control defined in \eqref{defhata}.  We can then estimate the value function at any time by regression similarly as in \eqref{regressV}.  Alternately, by considering the value function in the BSDE as in \eqref{BSDEY}, we can obtain an approximation of $V$  via the mean-field neural network $\Yc_i^*$ at any time $t_i$. 
Once again, in a backward resolution process, if the support of the distribution is unknown, it is necessary to employ an iterative algorithm, as suggested in Section \ref{subsec:MFNN}, to explore and identify the distributions of interest.

\subsection{Global algorithms}


In the spirit of the deep BSDE method in \cite{han2018solving}, we consider a mean-field neural network $\Uc_\theta(\mu)(x)$, and time dependent mean-field neural network $\Zc_\theta(t,\mu)(x)$, for approximating respectively the initial value of the $\Yc$ component, and the $\Zc$ component at any time of the MKV BSDE. 
We then define by forward induction: starting from $X_0$ $\sim$ $\mu_0$, $\Yc_0$ $=$ $\Uc_\theta(\mu_0)(X_0)$,  for $i$ $=$ $0,\ldots,N_T-1$, 
\begin{align}
X_{i+1} &= \;  X_i  + B(X_i,\mu_i,\Yc_i) \Delta t + \sigma(X_i,\mu_i) \Delta W_i,  \nonumber \\
\Yc_{i+1} &= \;  \Yc_i + \tilde\E\big[\Hc(X_i,\mu_i,\Yc_i,\Zc_\theta(t_i,\mu_i)(X_i),\tilde X_i,\tilde\Yc_i,\Zc_\theta(t_i,\mu_i)(\tilde X_i))\big] \Delta t + \Zc_\theta(t_i,\mu_i)(X_i)  \Delta W_i, \label{eq:forMKV}
\end{align}
and minimize over $\theta$ the global loss function
\beqs
L(\theta) &=& \E \Big| \Yc_{N_T} - G(X_{N_T},\mu_{N_T}) \Big|^2. 
\enqs
In practical implementation, we use a training sample of probability measures $\mu_0^{(m)}$, and then for each $m$, $N$ samples $X_0^{(m),(n)}$ $\sim$ $\mu_0^{(m)}$, $\Yc_0^{(m),(n)}$ $=$ $\Uc_0(\mu_0^{(m)})(X_0^{(m),(n)})$, $n$ $=$ $1,\ldots,N$, and  for $i$ $=$ $0,\ldots,N_T-1$
\begin{align}
X_{i+1}^{(m),(n)} &= \;  X_i^{(m),(n)}  + B(X_i^{(m),(n)},\hat\mu_i^{(m)},\Yc_i^{(m),(n)}) \Delta t + \sigma(X_i^{(m),(n)},\hat\mu_i^{(m),(n)}) \Delta W_i, \\
\Yc_{i+1}^{(m),(n)} &= \;  \Yc_i^{(m),(n)} 
+ \frac{\Delta t}{N}  \sum_{n'=1}^N  \Hc(X_i^{(m),(n)},\hat\mu_i^{(m)},\Yc_i^{(m),(n)},\Zc_\theta(t_i,\hat\mu_i^{(m)})(X_i^{(m),(n)}), \\
& \quad \hspace{1cm}  \tilde X_i^{(m),(n')},\tilde\Yc_i^{(m),(n')},\Zc_\theta(t_i,\hat\mu_i^{(m)})(\tilde X_i^{(m),(n')}))  \; + \;  \Zc_\theta(t_i,\hat\mu_i^{(m)})(X_i^{(m),(n)})  \Delta W_i,
\end{align}
where $\tilde X_i^{(m),(n)}$, $\tilde\Yc_i^{(m),(n)}$ are independent copies of $X_i^{(m),(n)}$, $\Yc_i^{(m),(n)}$, while $\hat\mu_0^{(m)}$ $=$ $\mu_0^{(m)}$, $\hat\mu_i^{(m)}$, $i$ $=$ $1,\ldots,N_T$, are estimated as in Section \ref{sec:globalcon}. 
The loss function is then approximated by 
\beqs
L_{M,N}(\theta) &=& \frac{1}{MN} \sum_{m=1}^M\sum_{n=1}^N \Big| \Yc_{N_T}^{(m),(n)} - G(X_{N_T}^{(m),(n)},\hat\mu_{N_T}^{(m)})  \Big|^2.    
\enqs

\vspace{1mm}

The pseudo-code is described in Algorithm \ref{algo:BSDEglobal}. 

\vspace{1mm}

\begin{algorithm}[H] 
{\bf Input data:} A mean-field neural network $\Uc_\theta(\mu)(x)$, and a time-dependent mean-field neural network $\Zc_\theta(t,\mu)(x)$.  

{\bf Initialization:} learning rate $\gamma$ and parameters $\theta$ 

\For{each epoch} 
{
Generate a batch  of $M$ initial distributions $\mu_0^{(m)}$, $m$ $=$ $1,\ldots,M$.\;
\For{each batch $m$}{
Generate Brownian increments $\Delta W_i^{(m),(n)}$, $i$ $=$ $0,\ldots,N_T-1$,  $n$ $=$ $1,\ldots,N$. \;  
Compute sample trajectories $X_i^{(m),(n)}$, $\tilde X_i^{(m),(n)}$, $\Yc_i^{(m),(n)}$, $\tilde\Yc_i^{(m),(n)}$, $n$ $=$ $1,\ldots,N$, and estimate $\hat\mu_i^{(m)}$, $i$ $=$ $0,\ldots,N_T$,  
}
Compute the batch loss $L_{M,N}(\theta)$ and its gradient $\nabla_\theta L_{M,N}(\theta)$ \;
Update $\theta$ $\leftarrow$ $\theta$ $-$ $\gamma \nabla_\theta L_{M,N}(\theta)$ \;
}
{\bf Return:} the set of optimized parameters $\theta^*$. 
\caption{Deep MKV BSDE} \label{algo:BSDEglobal} 
\end{algorithm}
\vspace{1mm}

The output of this global deep BSDE algorithm is an approximation of the $\Yc$ component of the BSDE at initial time $t_0$ $=$ $0$ by a mean-field neural network  $\Uc_{\theta^*}$, and yields approximation of the $\Zc$ component  at times $t_i$, $i$ $=$ $0,\ldots,N_T-1$, 
by mean-field neural networks $\Zc_{\theta^*}(t_i,\mu_i)(X_i)$ along the law  $\mu_i$, and  state $X_i$ explored during the learning algorithm.  
The value function can then be estimated at any time $t_k$  by regression as follows: we approximate the value function at time $t_k$ by a mean-field neural neural network $\vartheta_{\eta_k}(\mu)(x)$ valued in $\R$, and minimize over the parameters $\eta_k$ of this neural network the quadratic loss function 
 \begin{align}  \label{regressBSDEGlob} 
 \E   \Big|  Y_k  - \vartheta_{\eta_k}(\mu_k)(X_k)  \Big|^2, 
 \end{align} 
 where  $(X_k, Y_k)$ are generated by using equation \eqref{eq:forMKV} for $i=0, \ldots, k-1$ (here $Y_k$ is the first component of $\Yc_k$ $=$ $(Y_k,P_k)$ in \eqref{BSDEY}-\eqref{BSDEP}), 
 and $\mu_k$ is estimated from the distribution of the $X_k$.


\vspace{5mm}

In order to avoid the cost of solving equation \eqref{regressBSDEGlob} at each time step, we can propose two other global methods permitting to obtain directly the value function.

We first present a variation of the deep BSDE algorithm by considering two time-dependent mean-field neural networks $\Yc_\theta(t,\mu)(x)$ and $\Zc_\theta(t,\mu)(x)$, for approximating the pair solution of the MKV BSDE at any time. 
We then define by forward induction: starting from $X_0$ $\sim$ $\mu_0$,  for $i$ $=$ $0,\ldots,N_T-1$, 
\begin{align} \label{Xlocal}
X_{i+1} &= \;  X_i  + B(X_i,\mu_i,\Yc_\theta(t_i,\mu_i)(X_i)) \Delta t + \sigma(X_i,\mu_i) \Delta W_i, 
\end{align}
and minimize over $\theta$ the global loss function as a sum of local loss functions: 
\beqs
\tilde L(\theta) &=& \E\Big[ 
\sum_{i=1}^{N_T-1}   \Big| \Yc_{\theta}(t_{i+1},\mu_{i+1})(X_{i+1}) -  \Yc_{\theta}(t_i,\mu_i)(X_i) -  \Zc_{\theta}(t_i,\mu_i)(X_i) \Delta W_i \\
& & \quad \quad - \;   \tilde\E\big[\Hc(X_i,\mu_i,\Yc_{\theta}(t_i,\mu_i)(X_i),\Zc_{\theta}(t_i,\mu_i)(X_i),\tilde X_i,\Yc_{\theta}(t_i,\mu_i)(\tilde X_i),\Zc_{\theta}(t_i,\mu_i)(\tilde X_i))\big] \Delta t  \Big|^2 \Big] ,
\enqs
with the convention that  $\Yc_\theta(t_{N_T},\mu)(x)$ $=$ $G(x,\mu)$.  In practical implementation, we use a training sample of probability measures $\mu_0^{(m)}$, and then for each $m$ $=$ $1,\ldots,M$,  $N$ samples $X_0^{(m),(n)}$ $\sim$ $\mu_0^{(m)}$, $X_i^{(m),(n)}$, $\tilde X_i^{(m),(n)}$, $n$ $=$ $1,\ldots,N$, according to \eqref{Xlocal}, 
and estimated probability measures $\hat\mu_i^{(m)}$, $i$ $=$ $1,\ldots,N_T$.  The loss function is then approximated by $\tilde L_{M,N}(\theta)$.

\vspace{1mm}

The pseudo-code is described in Algorithm \ref{algo:BSDEglobal2}. 

\vspace{1mm}

\begin{algorithm}[H] 
{\bf Input data:}  Two time-dependent mean-field neural network $\Yc(t,\mu)(x)$, $\Zc_\theta(t,\mu)(x)$.  

{\bf Initialization:} learning rate $\gamma$ and parameters $\theta$ \;
\For{each epoch} 
{
Generate a batch  of $M$ initial distributions $\mu_0^{(m)}$, $m$ $=$ $1,\ldots,M$  \; 
\For{each batch $m$}
{
Generate Brownian increments $\Delta W_i^{(m),(n)}$, $i$ $=$ $0,\ldots,N_T-1$,  $n$ $=$ $1,\ldots,N$. \;  
Compute sample trajectories $X_i^{(m),(n)}$, $\tilde X_i^{(m),(n)}$,  $n$ $=$ $1,\ldots,N$, and estimate $\hat\mu_i^{(m)}$, $i$ $=$ $0,\ldots,N_T$,  
}
Compute the batch loss $\tilde L_{M,N}(\theta)$ and its gradient $\nabla_\theta \tilde L_{M,N}(\theta)$ \;
Update $\theta$ $\leftarrow$ $\theta$ $-$ $\gamma \nabla_\theta\tilde L_{M,N}(\theta)$ \;
}
{\bf Return:} the set of optimized parameters $\theta^*$. 
\caption{Deep MKV BSDE global/local } \label{algo:BSDEglobal2} 
\end{algorithm}


\vspace{5mm}

Finally, we present a multi-step version of the deep MKV BSDE algorithm.  We consider two time-dependent mean-field neural networks $\Yc_\theta(t,\mu)(x)$ and $\Zc_\theta(t,\mu)(x)$, for approximating the pair solution of the MKV BSDE at any time, and  
define by forward induction: starting from $X_0$ $\sim$ $\mu_0$,  for $i$ $=$ $0,\ldots,N_T-1$, 
\begin{align} \label{Xmulti} 
X_{i+1} &= \;  X_i  + B(X_i,\mu_i,\Yc_\theta(t_i,\mu_i)(X_i)) \Delta t + \sigma(X_i,\mu_i) \Delta W_i. 
\end{align}
The global  loss function to be minimized is of the form
 
 \beqs
L^{multi}(\theta) &=& \E \Big[ \sum_{i=0}^{N_T-1} \Big| G(X_{N_T},\mu_{N_T}) -   \sum_{j=i}^{N_T-1}  \Zc_\theta(t_j,\mu_j)(X_j) \Delta W_j  -  \Yc_{\theta}(t_i,\mu_i)(X_i)    \\
 & &  -  \; \sum_{j=i}^{N_T-1}  \tilde\E\big[\Hc(X_j,\mu_j,\Yc_{\theta}(t_j,\mu_j)(X_j),\Zc_{\theta}(t_j,\mu_j)(X_j),\tilde X_j,\Yc_{\theta}(t_j,\mu_j)(\tilde X_j),\Zc_{\theta}(t_j,\mu_j)(\tilde X_j))\big] \Delta t 
 \Big|^2 \Big].  
 \enqs
 
Again,  in practical implementation, we use a training sample of probability measures $\mu_0^{(m)}$, and then for each $m$ $\in$ $\{1,\ldots,M\}$,  $N$ samples $X_0^{(m),(n)}$ $\sim$ $\mu_0^{(m)}$, $X_i^{(m),(n)}$, $\tilde X_i^{(m),(n)}$, $n$ $=$ $1,\ldots,N$, according to \eqref{Xmulti},  
and estimated probability measures $\hat\mu_i^{(m)}$, $i$ $=$ $1,\ldots,N_T$.  The loss function is then approximated by $L^{multi}_{M,N}(\theta)$.

\vspace{2mm}

The pseudo-code is described in Algorithm \ref{algo:BSDEglobalmulti}.

\vspace{1mm}

\begin{algorithm}[H] 
{\bf Input data:}  Two time-dependent mean-field neural networks $\Yc(t,\mu)(x)$, $\Zc_\theta(t,\mu)(x)$.  

{\bf Initialization:} learning rate $\gamma$ and parameters $\theta$ \;
\For{each epoch} 
{
Generate a batch  of $M$ initial distributions $\mu_0^{(m)}$, $m$ $=$ $1,\ldots,M$  \; 
\For{each batch $m$}
{
Generate Brownian increments $\Delta W_i^{(m),(n)}$, $i$ $=$ $0,\ldots,N_T-1$,  $n$ $=$ $1,\ldots,N$. \;  
Compute sample trajectories $X_i^{(m),(n)}$, $\tilde X_i^{(m),(n)}$,  $n$ $=$ $1,\ldots,N$, and estimate $\hat\mu_i^{(m)}$, $i$ $=$ $0,\ldots,N_T$,  
}
Compute the batch loss $L^{multi}_{M,N}(\theta)$ and its gradient $\nabla_\theta L^{multi}_{M,N}(\theta)$ \;
Update $\theta$ $\leftarrow$ $\theta$ $-$ $\gamma \nabla_\theta L^{multi}_{M,N}(\theta)$ \;
}
{\bf Return:} the set of optimized parameters $\theta^*$. 
\caption{Deep multi-step MKV BSDE  } \label{algo:BSDEglobalmulti} 
\end{algorithm}

\vspace{2mm}

The output of  Algorithms \ref{algo:BSDEglobal2} and \ref{algo:BSDEglobalmulti} is an approximation of the $\Yc$ component of the BSDE at initial time $t_0$ $=$ $0$ by a mean-field neural network  $\Yc_{\theta^*}(t_0,.)(.)$, and yields approximation of the $\Yc$, at other times $t_i$, $i$ $=$ $1,\ldots,N_T-1$, and 
$\Zc$ at times  $t_i$, $i$ $=$ $0,\ldots,N_T-1$, by mean-field neural networks $\Yc_{\theta_i^*}(t_i,\mu_i)(X_i)$,  $\Zc_{\theta^*}(t_i,\mu_i)(X_i)$ along the law  $\mu_i$, and  state $X_i$ explored during the learning algorithm.

In the case of global algorithms using the cylindrical network, there is no requirement to adapt the parameter $\Kc$. The need for adaptation methods, as proposed in Section \ref{subsec:MFNN}, arises primarily when employing the bin method.

\section{Numerical examples} \label{secnum}

We shall illustrate the results of our different algorithms on three test cases.  The two first examples are MKV control problems where the diffusion coefficient is constant, and the BSDE approach can be used.  The third example is  
a classical mean variance problem, hence with control on the diffusion coefficient.  
We then test the three cases using the dynamic programming-based algorithms and  for  the two first cases using also  the backward SDE-based algorithms.

For each problem, we will test the optimized solutions $v(\mu_0)$ found by using different initial distributions $\mu_0$ and compare the result obtained to the analytical solution or the reference calculated by an other method. 
For all test cases, we keep the same parameters for the neural networks:
\begin{itemize}
\item For the bin method, we take 2 layers of 20 neurons.
\item For the cylinder method, we take 2 layers of 20 neurons for the two networks.
\end{itemize}
For both methods we use the $\tanh$ activation function. At each iteration of the ADAM gradient method \cite{kingma2014adam}, we consider for each of the $M$ tested distributions  $N=100000$ realizations of the process $X$. 
These parameters are chosen accordingly the results of \cite{phawar22a}.  We  either take  a batch size equal to $M=5$, $M=8$, $M=10$ or $M=20$, using between $30000$ to $120000$ gradient iterations: we have to adapt the batch size and the number of gradient iterations to be able to solve the problem on the  graphic card GPU NVidia V100 32Gb (except when specified due to  memory limitation) and in order to obtain the result in less than 3 days. $K$ in the tables below is the number of bins used, and $\Delta t$ $=$ $T/N_T$ is the time step.

\subsection{The test examples}

\subsubsection{Systemic risk model}
\label{sec:Systemic}

We consider a mean-field model of systemic risk introduced in \cite{carfousun}.  This model was introduced in the context of mean field games but here we consider a cooperative version. The limit problem (when the number of banks is large) of the social planner (central bank) is formulated as follows. The log-monetary reserve of the representative bank is governed by the mean-reverting controlled McKean-Vlasov dynamics 
\begin{align}
dX_t &= \; \big[  \kappa( \E[X_t] - X_t) + \alpha_t]  \ \di t +  \sigma d W_t , 
\quad X_0 \sim \mu_0, 
\end{align}
where $\alpha$ $=$ $(\alpha_t)_t$  is the  control rate of borrowing/lending to a central bank that aims to minimize  the functional cost
\begin{align} \label{defV0}
J(\alpha) \; = \; \E \Big[ \int_0^T \tilde f(X_t,\E[X_t],\alpha_t) \ \di t + \tilde g(X_T,\E[X_T]) \Big] &\;  \rightarrow \quad v(\mu_0)  \; = \; \inf_{\alpha} J(\alpha), 
\end{align}
where the running and terminal costs are given by 
\begin{align}
\tilde f(x,\bar x,a)  \; = \; \frac{1}{2}a^2 - qa(\bar x-  x) + \frac{\eta}{2} (\bar x-x)^2, & \quad \tilde g(x,\bar x) \; = \; \frac{c}{2}(x-\bar x)^2, 
\end{align} 
for some positive constants $q$, $\eta$, $c$ $>$ $0$, with $q^2$ $\leq$ $\eta$.  
Notice that in this linear-quadratic example, the objective function is convex with respect to the control process, which ensures  the  convergence of the global algorithm.

The explicit solution of the linear-quadratic McKean-Vlasov control problem \eqref{defV0} is solved via the resolution of a Riccati equation (see \cite{baspha19}), and is analytically given by 
\begin{align}
v(t,\mu) &=  \int_\R V(t, x, \mu) \mu(\d x) \; =  \; Q_t  \int_\R (x - \bar \mu)^2 \mu(\d x) 
+ \sigma^2 \int_t^T Q_s \ \di s,
\label{eq:vAnal}
\end{align}
where we set $\bar\mu$ $:=$ $\E_{\xi\sim\mu}[\xi]$ $=$ $\int_\R x \mu(\d x)$, and 
\begin{align}
Q_t & = 
- \frac{1}{2} \Big[ \kappa + q - \sqrt{\Delta} \frac{ \sqrt{\Delta} \sinh(\sqrt{\Delta}(T-t))  + (\kappa + q + c) \cosh (\sqrt{\Delta}(T-t))}{ \sqrt{\Delta} \cosh(\sqrt{\Delta}(T-t))  + 
(\kappa + q + c) \sinh (\sqrt{\Delta}(T-t))} \Big],   
\end{align}     
with $\sqrt{\Delta}$ $=$ $\sqrt{(\kappa+q)^2  + \eta -  q^2}$,  and 
\begin{align}
\int_t^T Q_s \ \di s & = \;  \frac{1}{2} \ln \Big[ \cosh (\sqrt{\Delta}(T-t)) + \frac{\kappa + q  + c}{\sqrt{\Delta}}  \sinh (\sqrt{\Delta}(T-t)) \Big] - \frac{1}{2} (\kappa + q)(T-t). 
\end{align} 
In this example, the function $\hat\mra$ that attains the infimum of the Hamiltonian function is $\hat\mra(x,\mu,p)$ $=$ $q(\bar\mu -x) - p$,  the function in \eqref{defUc} is $\Uc(t,x,\mu)$ $=$ $2Q_t(x-\bar\mu)$, which yields the optimal feedback control: 
$\mfa^\star(t,x,\mu)$ $=$ $(q+2Q_t)(\bar\mu - x)$.  The  BSDE  \eqref{BSDEY}-\eqref{BSDEP}  is then written as 
\begin{equation}
\left\{ 
\begin{array}{ccl}
\d X_t &=&   \big[ (\kappa +q) (\E[X_t ]- X_t)  - P_t ] \d t + \sigma \d W_t, \quad X_0 \sim \mu_0, \\
\d Y_t & = & - \big[   \frac{1}{2}( \eta -q^2) (\E[X_t ]- X_t)^2 +  \frac{1}{2} P_t^2 \big] \d t + Z_t \d W_t, \quad Y_T = \frac{c}{2}(X_T - \E[X_T])^2,   \\
\d P_t & = & \big[ - (\kappa + q)  (\E[P_t ]- P_t)  +(\eta - q^2) (\E[X_t ]- X_t) \big] \d t + M_t \d W_t, \quad  P_T = - c( \E[X_T] - X_T). 
\end{array}
\right.
\end{equation}

\vspace{3mm}

\begin{figure}[H]
\begin{minipage}[t]{0.30\linewidth}
  \centering
 \includegraphics[width=\textwidth]{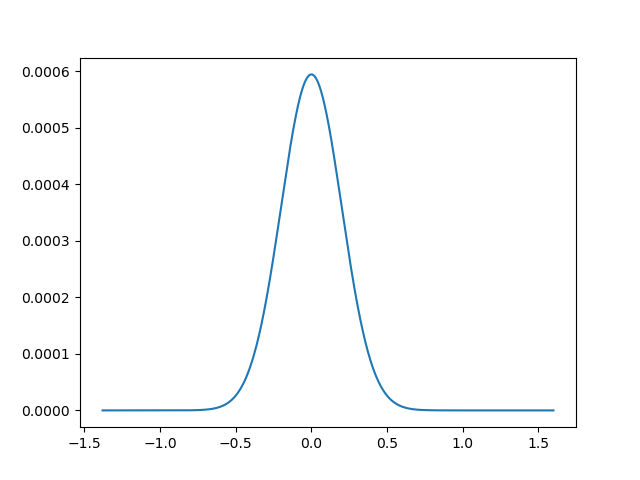}
\caption*{Case 1}
\end{minipage}
\begin{minipage}[t]{0.30\linewidth}
  \centering
 \includegraphics[width=\textwidth]{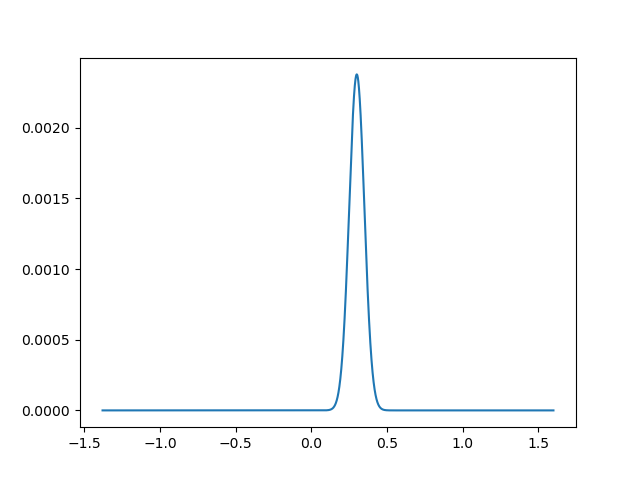}
\caption*{ Case 2}
\end{minipage}
\begin{minipage}[t]{0.30\linewidth}
  \centering
 \includegraphics[width=\textwidth]{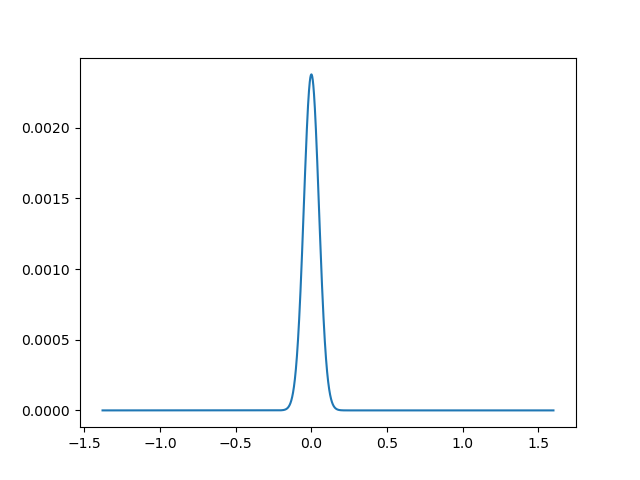}
\caption*{Case 3 }
\end{minipage}
\begin{minipage}[t]{0.30\linewidth}
  \centering
 \includegraphics[width=\textwidth]{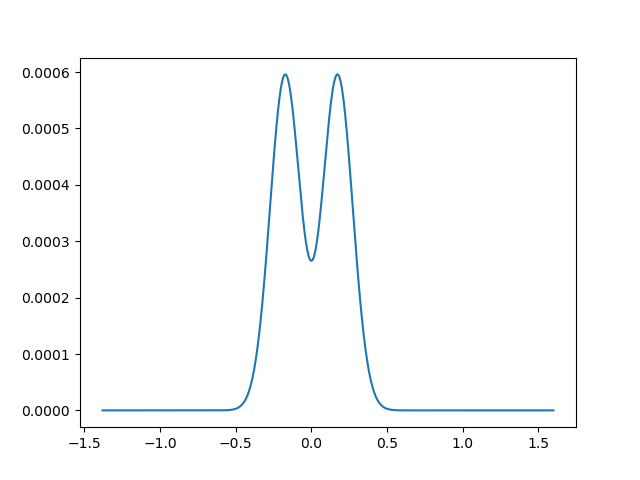}
\caption*{ Case 4}
\end{minipage}
\begin{minipage}[t]{0.30\linewidth}
  \centering
 \includegraphics[width=\textwidth]{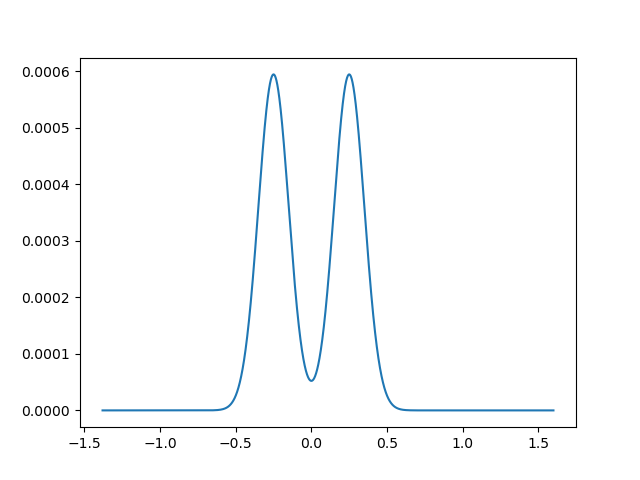}
\caption*{ Case 5}
\end{minipage}
\begin{minipage}[t]{0.30\linewidth}
  \centering
 \includegraphics[width=\textwidth]{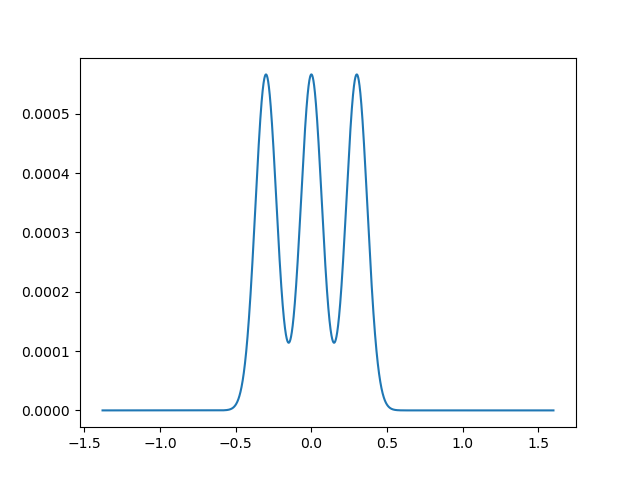}
\caption*{ Case 6}
\end{minipage}
\caption{Distribution $\mu_0$ tested on the systemic case. \label{fig:systemicDist} }
\end{figure}

For the numerical tests of  the different  methods, we take   $\sigma=1$, $\kappa =0.6$, $q=0.8$,  $T=0.2$, $C=2$, $\eta=2$. We solve the problem \eqref{defV0} using our various algorithms and compare  the solution  obtained at $t$ $=$ $0$ with  $v(0,\mu_0)$  given by \eqref{eq:vAnal} for different initial distributions 
$\mu_0$ plotted on Figure \ref{fig:systemicDist}:
\begin{itemize}
    \item Case 1 : Gaussian with $ \bar \mu_0 =0$, $std(\mu_0)=0.2$,
    \item Case 2 : Gaussian  with $ \bar \mu_0 =0.3$, $std(\mu_0)=0.05$,
    \item Case 3 : Gaussian with $ \bar \mu_0 =0.$, $std(\mu_0)=0.05$,
    \item Case 4 : Mixture of two Gaussian random variables: $X_0 =  P (-k + \theta Y) + (1-P) (k + \theta \bar Y)$ with $P$ a  Bernouilli random variable  with parameter $\frac{1}{2}$, $k=\frac{\sqrt{3}}{10}$, $\theta = 0.1$, $Y, \bar Y \sim  \mathcal{N}(0,1)$,
    \item Case 5 : Mixture of two Gaussian random variables $X_0 =  P (-k + \theta Y) + (1-P) (-k + \theta \bar Y)$ with $P$ a Bernouilli random variable with parameter $\frac{1}{2}$, $k=0.25$, $\theta = 0.1$, $Y, \bar Y \sim  \mathcal{N}(0,1)$,
    \item Case 6 : Mixture of 3 Gaussian random variables : $X_0 =  [ - 1_{\lfloor 3U \rfloor= 0}  k + 
1_{\lfloor 3U \rfloor= 1 } k ]+  \theta Y $ with $U \sim U(0,1)$, $k=0.3$, $\theta = 0.07$, $\bar Y \sim  \mathcal{N}(0,1)$.  
\end{itemize}

Notice that case 1 and 4 have the same variance for $\mu_0$ so that the values $v(0,\mu_0)$ of \eqref{eq:vAnal} should be the same. Similarly, values of case 2 and 3 are the same.

\subsubsection{Min/max linear quadratic MKV control}
\label{sec:MinMax}
 
We next  consider a mean-field model in which the dynamics is linear, the running cost is quadratic in the position, the control and the expectation of the position, while  the terminal cost gives inventive to be close to one of two targets. This type of model is inspired by the min-LQG problem of~\cite{salhab2015a}. More precisely, we consider the following controlled McKean-Vlasov dynamics 
\begin{align}
    \di X_t &= \; \big[  A X_t + \bar A \E[X_t] + B \alpha_t \big]  \ \di t +  \sigma \ \di W_t , 
\quad X_0 \sim \mu_0, 
\end{align}
where $\alpha$ $=$ $(\alpha_t)_t$  is the  control, and the agent aims to minimize  the functional cost
\begin{align*} 
 \label{defV0-minLQ}
    J(\alpha) \; = \; \E \Big[ \int_0^T f(X_t,\E[X_t],\alpha_t) \ \di t + g(X_T) \Big] 
    &\;
    \rightarrow 
    \quad v(\mu_0) \; = \; \inf_{\alpha} J(\alpha), 
\end{align*}
where the running and terminal costs are given by 
\begin{align}
    f(x,\bar x,a)  \; = \; \frac{1}{2} \left( Q x^2 + \bar Q (x - S \bar x)^2 + R a^2 \right), 
    & \quad g(x) \; = \; \min \left\{  |x - \zeta_1 |^2, | x - \zeta_2 |^2 \right\}, 
\end{align} 
for some non-negative constants $Q$, $\bar Q$, $S$, $R$, and two real numbers $\zeta_1$ and $\zeta_2$.  
Notice that $g$ is not a convex function, and the solution to the MKV BSDE is not necessarily an optimal control.

\vspace{1mm}

In this example, the BSDE  \eqref{BSDEY}-\eqref{BSDEP}  is then written as 
\begin{equation}
\left\{
\begin{array}{ccl}
\d X_t &= & [AX_t + \bar A\E[X_t]  - \frac{B^2}{R} P_t ] \d t + \sigma \d W_t, \quad X_0 \sim \mu_0 \\
\d Y_t & = & -  \frac{1}{2}\big[ Q X_t^2  + \bar Q(X_t - S \E[X_t])^2 + \frac{B^2}{R} P_t^2 \big] \d t + Z_t \d W_t, \quad Y_T =  \min[ |X_T - \zeta_1|^2, |X_T - \zeta_2|^2 ] \\
\d P_t & = & - \big[ AP_t + \bar A \E[P_t] +  Q X_t + \bar Q(X_t - \E[X_t]) + \bar Q (S-1)^2\E[X_t] \big] \d t + M_t \d W_t,  \\
P_T &= &  2 \big(X_T - \min(\zeta_1,\zeta_2) 1_{X_T \leq \frac{\zeta_1+\zeta_2}{2}} -   \max(\zeta_1,\zeta_2) 1_{X_T > \frac{\zeta_1+\zeta_2}{2}} \big). 
\end{array}
\right.
\end{equation}

\vspace{1mm}

\begin{figure}[H]
\begin{minipage}[t]{0.32\linewidth}
  \centering
 \includegraphics[width=\textwidth]{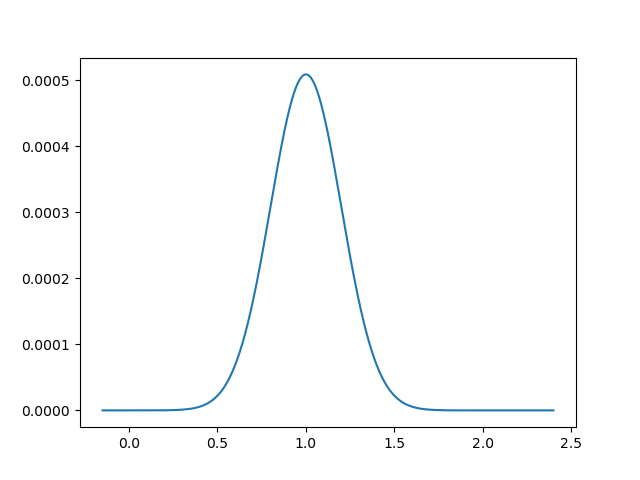}
\caption*{ Case 1}
\end{minipage}
\begin{minipage}[t]{0.32\linewidth}
  \centering
 \includegraphics[width=\textwidth]{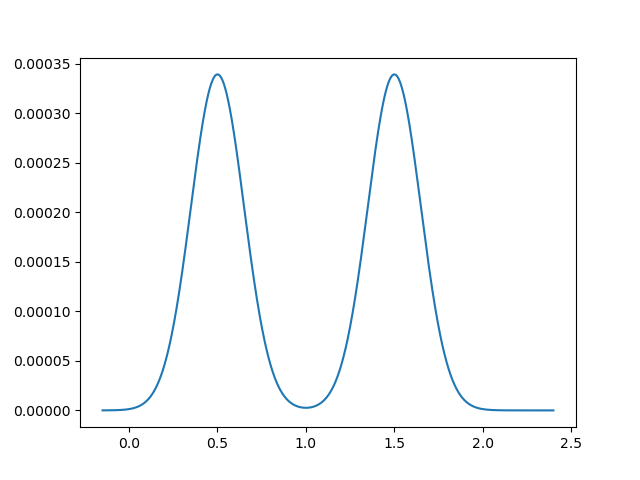}
\caption*{ Case 2}
\end{minipage}
\begin{minipage}[t]{0.32\linewidth}
  \centering
 \includegraphics[width=\textwidth]{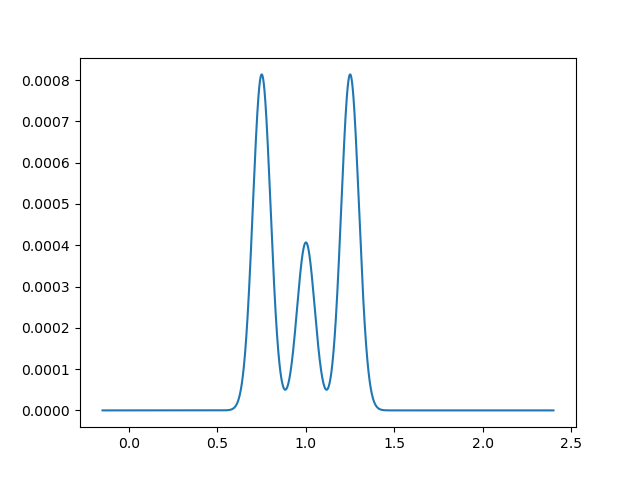}
\caption*{ Case 3 }
\end{minipage}
\caption{Distribution $\mu_0$ tested on the min/max linear case. \label{fig:distMinDist}}
\end{figure}

For the numerical  tests, we take $A= 1$, $\bar A = 0.5$, $B=1$, $Q= \bar Q = R =S =1$, $\sigma = 0.5$, $\zeta_1=0.25$, $\zeta_2=1.75$. 
We first solve the problem \eqref{MKVcontrol}  by the different algorithms and we can compare the solution $v(\mu_0)$  obtained  for different distributions $\mu_0$ to a reference  calculated using \cite{carlau19} approach. Notice that \cite{carlau19} method needs to be run for each initial distribution tested. 
We use three different distributions $\mu_0$ plotted on Figure \ref{fig:distMinDist}:
\begin{itemize}
    \item Case 1 : Gaussian distribution $ \bar \mu_0 =1$, $std(\mu_0)=0.2$. The reference values are   $0.484$ for $T=0.2$, and $0.818$ for $T=0.5$. 
    \item Case 2 : Mixture of  two Gaussian random variables : $X_0=  P ( \zeta_1 +\theta Y) + (1-P) ( \zeta_2 +  \theta \bar Y)$ with $P$ a Bernouilli random variable with  parameter $\frac{1}{2}$, $\theta = 0.15$, $Y , \bar Y, \tilde Y \sim  \mathcal{N}(0,1)$, with reference values   $0.494$ for $T=0.2$, and  $1.082$ for $T=0.5$.
    \item Case 3 : Mixture of three Gaussian random variables: $X_0 = [  1_{\lfloor 5U \rfloor <  2}  \zeta_1 + 
1_{\lfloor 5U \rfloor >  3 } \zeta_2 + 1_{ 2 \le \lfloor 5U \rfloor \le  3} (\zeta_1+\zeta_2) ]+  \theta Y $ with $U \sim U(0,1)$, $\theta = 0.05$ with reference values $ 0.491$  for $T=0.2$, and $0.836$ for $T=0.5$.
\end{itemize}

\subsubsection{Mean-variance problem} 
\label{sec:MeanVar}

We consider the celebrated Markowitz portfolio selection problem where an investor can invest at any time $t$ an amount $\alpha_t$ in a risky asset (assumed for simplicity to follow a Black-Scholes model with constant rate of return $\beta$ and volatility $\nu$ $>$ $0$), hence generating a wealth process $X$ $=$ $X^\alpha$ with dynamics
\begin{align*}
dX_t &= \; \alpha_t \beta dt + \alpha_t \nu dW_t, \quad 0 \leq t \leq T, \quad X_0 \sim \mu_0.   
\end{align*}
The goal is then to minimize over portfolio control $\alpha$ the mean-variance criterion: 
\begin{align*}
J(\alpha) &= \; \lambda {\rm Var}(X_T^\alpha) - \E[X_T^\alpha],
\end{align*}
 where $\lambda$ $>$ $0$ is a parameter related to the risk aversion of the investor.  
 
We refer to \cite{ismpha19} for the McKean-Vlasov approach to Markowitz mean-variance problems (in a more general context), and we recall that the solution to the Bellman equation is given by  
\begin{align} \label{eq:analMV}
V(t,x,\mu) & = \; \lambda e^{-R(T-t)} (x - \bar\mu)^2 - x - \frac{1}{4\lambda} \big[ e^{R(T-t)} - 1\big], \\
\Uc(t,x,\mu) &  = \; 2 \lambda e^{-R(T-t)}  (x - \E_\mu[\xi]) - 1,
\end{align}
where we set $R$ $:=$ $\beta^2/\nu^2$. 
Moreover,  the optimal feedback control is given by 
\begin{align}
\mfa^*(t,x,\mu) &= \; - \frac{\beta}{\nu^2} \big(x - \bar\mu - \frac{e^{R(T-t)}}{2\lambda} \big). 
\label{eq:optCont}
\end{align}
Note that with this model, the BSDE approach cannot be used as the  volatility is controlled.

\begin{figure}[H]
\begin{minipage}[t]{0.30\linewidth}
  \centering
 \includegraphics[width=\textwidth]{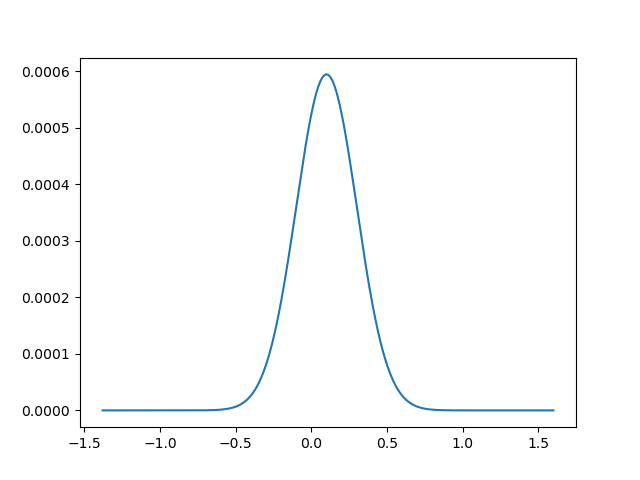}
\caption*{  Case 1}
\end{minipage}
\begin{minipage}[t]{0.30\linewidth}
  \centering
 \includegraphics[width=\textwidth]{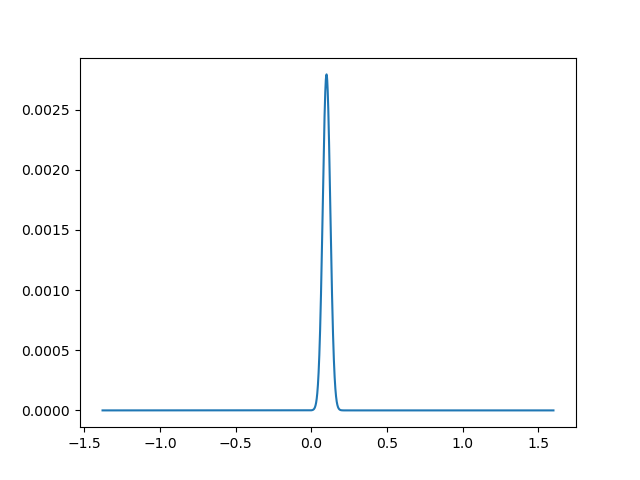}
\caption*{ Case 2}
\end{minipage}
\begin{minipage}[t]{0.30\linewidth}
  \centering
 \includegraphics[width=\textwidth]{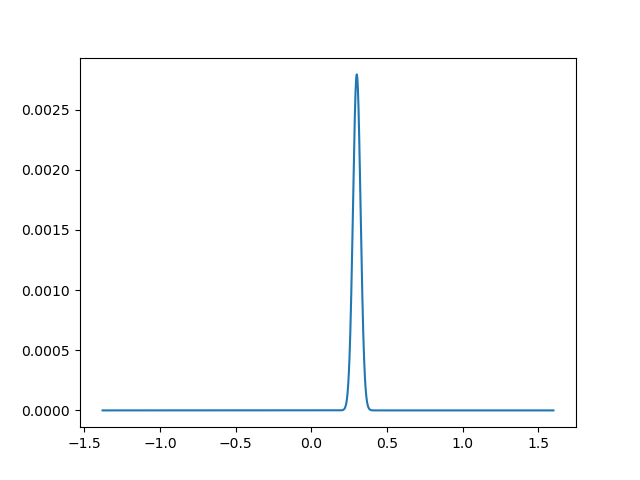}
\caption*{  Case 3 }
\end{minipage}
\begin{minipage}[t]{0.30\linewidth}
  \centering
 \includegraphics[width=\textwidth]{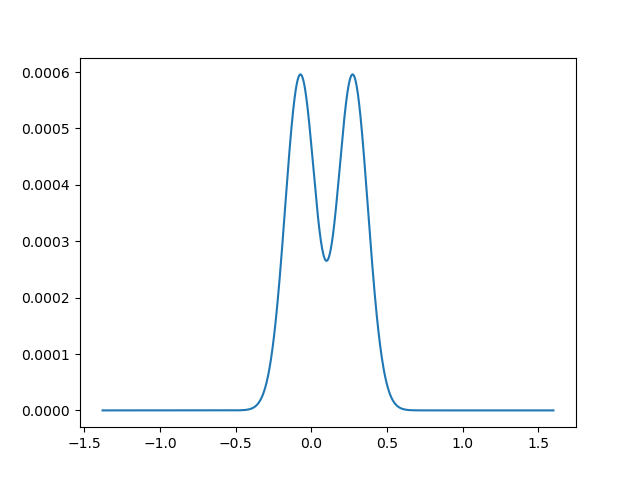}
\caption*{  Case 4}
\end{minipage}
\begin{minipage}[t]{0.30\linewidth}
  \centering
 \includegraphics[width=\textwidth]{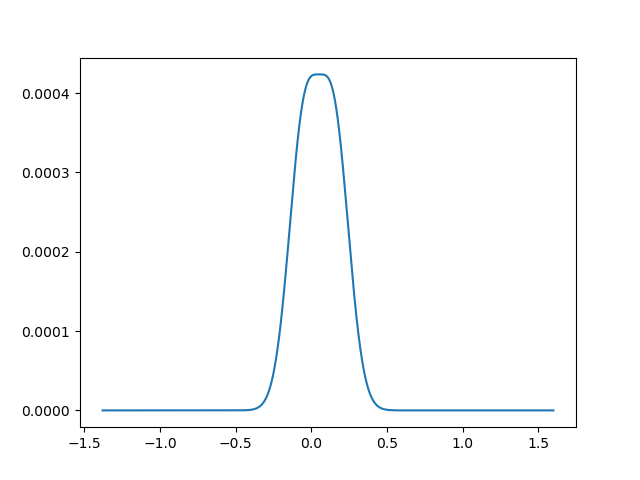}
\caption*{  Case 5}
\end{minipage}
\begin{minipage}[t]{0.30\linewidth}
  \centering
 \includegraphics[width=\textwidth]{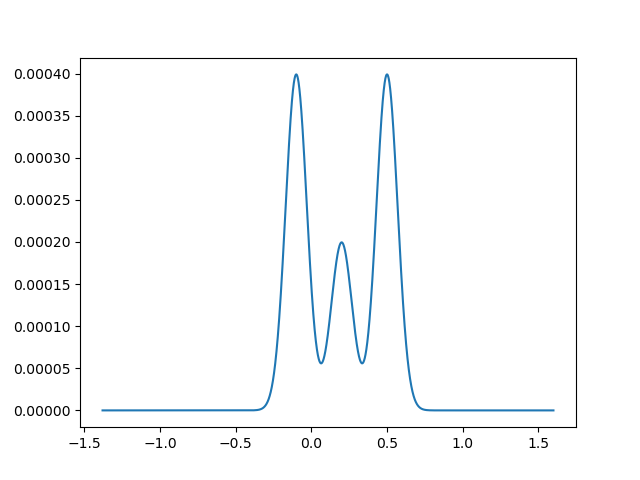}
\caption*{  Case 6}
\end{minipage}
\caption{Distribution $\mu_0$ tested on the mean variance case. \label{fig:distMeanVar}}
\end{figure}

We test our algorithms with the parameters $\beta = 0.1$, $ \nu = 0.4$, $\lambda = 0.5$.
We compare the solutions obtained at $t$ $=$ $0$ to the analytical solution $v(\mu_0)$ $=$ $\E_{\xi\sim\mu_0}[V(0,\xi,\mu_0)]$ given by \eqref{eq:analMV} for different initial distributions $\mu_0$ plotted in Figure \ref{fig:distMeanVar}, and explicitly given by: 
\begin{itemize}
    \item Case 1 : Gaussian  distribution with $ \bar \mu_0 =0.1$, $std(\mu_0)=0.2$.
    \item Case 2 : Gaussian  distribution with $ \bar \mu_0 =0.2$, $std(\mu_0)=0.025$. 
    \item Case 3 : Gaussian distribution with $ \bar \mu_0 =0.3$, $std(\mu_0)=0.025$.
    \item Case 4 : Mixture of two Gaussian random variables: $X_0 =  P (-k +a + \theta Y) + (1-P) (-k + a+ \theta \bar Y)$ with $P$ a Bernouilli random variable with  parameter $\frac{1}{2}$, $k=\frac{\sqrt{3}}{10}$,$a=0.1$, $\theta = 0.1$, $Y , \bar Y \sim  \mathcal{N}(0,1)$,
    \item Case 5 : Mixture of two Gaussian random variables: $X_0 =  P (-k +a+ \theta Y) + (1-P) (-k + a+ \theta \bar Y)$ with $P$ a Bernouilli random variable with  parameter $\frac{1}{2}$, $a=0.05$, $k=0.1$, $\theta = 0.1$, $Y , \bar Y \sim  \mathcal{N}(0,1)$,
    \item Case 6 : Mixture of 3 Gaussian random variables: $X_0 =  a+ [ - 1_{\lfloor 5U \rfloor <  2}  k + 
1_{\lfloor 5U \rfloor >  3 } k ]+  \theta Y $ with $U \sim U(0,1)$, $a=0.2$, $k=0.3$, $\theta = 0.07$,$\bar Y \sim  \mathcal{N}(0,1) $. 
\end{itemize}

\subsubsection{A toy example of non LQ MKV control problem}

 We consider a one-dimensional controlled mean-field   dynamics of the form
 \beqs
\d X_t &=& \big[ \beta(X_t,\P_{X_t}) + \alpha_t \big]  \d t +   \sigma \d W_t,  \quad 0 \leq t\leq T, \quad X_0 \sim \mu_0,  
\enqs 
with  a cost functional of the form 
\begin{align*}
J(\alpha) &=\;  \E \Big[ \int_0^T  \big(  F(t,X_t,\P_{X_t}) + \frac{1}{2} |\alpha_t|^2  \big)  \d t + g(X_T,\P_{X_T}) \Big], \quad \rightarrow \quad v(\mu_0) \; = \; \inf_{\alpha\in \Ac} J(\alpha),  
\end{align*}   
where $g$ is  of the form:
\beqs
g(x,\mu) &=& \E_{\xi\sim\mu}[w(x-\xi)], 
\enqs
for some smooth $C^2$ even function $w$ on $\R$, e.g. $w(x)$ $=$ $\cos(x)$, and $F$ is a function to be chosen later.  

In this case, the  optimal feedback control   valued in $A$ $=$ $\R$ is given by
\beqs
\mfa^\star(t,x,\mu) &=& \hat\mra(t,x,\Uc(t,x,\mu)) \; = \;  - \Uc(t,x,\mu) \; = \; -  \partial_\mu v(t,\mu)(x)  \;\;  \mbox{ with } v(t,\mu) \; = \;  \E_{\xi\sim\mu}[V(t,\xi,\mu)],
\enqs
and  $V$ is solution to the Master Bellman equation: 
\begin{align}
\partial_t V(t,x,\mu) + \big( \beta(x,\mu) - \Uc(t,x,\mu) \big)  \partial_x V(t,x,\mu)  +   \frac{\sigma^2}{2}  \partial_{xx}^2 V(t,x,\mu) & \nonumber \\
+ \;  \E_{\xi\sim\mu} \Big[ \big( \beta(\xi,\mu) - \Uc(t,\xi,\mu) \big)  \partial_\mu V(t,x,\mu)(\xi)   +   \frac{\sigma^2}{2}  \partial_{x'}\partial_\mu V(t,x,\mu)(\xi)  \Big] & \nonumber \\
\quad + \;   F(t,x,\mu) + \frac{1}{2}   | \Uc(t,x,\mu)|^2   & = \; 0, \label{HJB} 
\end{align} 
with the terminal condition $V(T,x,\mu)$ $=$ $g(x,\mu)$. 

We look for a solution to the Master equation of the form: $V(t,x,\mu)$ $=$ $e^{T-t}\E_{\xi\sim\mu}[w(x-\xi)]$.  For such function $V$, we have $\partial_t V(t,x,\mu)$ $=$ $-V$, 
\beqs
\partial_x V(t,x,\mu) &=& e^{T-t} \E_{\xi\sim\mu} [ w'(x-\xi) ], \quad   \partial_{xx}^2 V(t,x,\mu) \; = \;  e^{T-t} \E_{\xi\sim\mu} [ w''(x-\xi) ] \\
\partial_\mu V(t,x,\mu)(\xi) &=& -  e^{T-t}  w'(x-\xi), \quad \partial_{x'}\partial_\mu V(t,x,\mu)(\xi) \; = \; e^{T-t}  w''(x-\xi), 
\enqs
and 
\beqs
\Uc(t,x,\mu) &=& e^{T-t} \E_{\xi\sim\mu} [ w'(x-\xi) - w'(\xi- x) ] \; = \; 2  e^{T-t} \E_{\xi\sim\mu} [ w'(x-\xi) ]  \; = \; 2 \partial_x V(t,x,\mu). 
\enqs
since $w$ is even. By plugging these derivatives expressions of $V$  into the l.h.s. of \eqref{HJB}, we then see that by choosing $F$ equal to 
\beqs
F(t,x,\mu) &=& e^{T-t} \E_{\xi\sim\mu} \Big[ (w- \sigma^2 w'')(x-\xi) + ( \beta(\xi,\mu) - \beta(x,\mu)) w'(x-\xi) \Big] \\
& & \; -  \;  2 e^{2(T-t)}  \E_{(\xi,\xi')\sim\mu\otimes\mu} \big[ w'(x-\xi)  w'(\xi- \xi') \big],
\enqs
the function $V$ satisfies the Master Bellman equation. 

For the choice of $w(x)$ $=$ $\cos(x)$, and using trigonometric relations, the function $F$ is written as
\beqs
F(t,x,\mu) &=& \cos(x)  [e^{T-t} \left( (1 + \sigma^2) \E_{\xi\sim\mu}( \cos(\xi)) + \E_{\xi\sim\mu}( \sin(\xi) \beta(\xi, \mu)) - \beta(x, \mu) \E_{\xi\sim\mu}(\sin(\xi))\right) - \\
& & 2 e^{2(T-t)} \left( \E_{\xi\sim\mu}( \sin(\xi) \cos(\xi)) \E_{\xi\sim\mu}(\sin(\xi))  -\E_{\xi\sim\mu}(\sin^2(\xi)) \E_{\xi\sim\mu}(\cos(\xi)) \right) ] + \\
&& \sin(x) [ e^{T-t}  \left(  (1 + \sigma^2)  \E_{\xi\sim\mu}(\sin(\xi)) -\E_{\xi\sim\mu}(\beta(\xi, \mu) \cos(\xi)) +\beta(x ,\mu) \E_{\xi\sim\mu}(\cos(\xi)) \right) - \\
& & 2e ^{2(T-t)} \left(\E_{\xi\sim\mu}(\sin(\xi)\cos(\xi)) \E_{\xi\sim\mu}(\cos(\xi) )-\E_{\xi\sim\mu}(\cos^2(\xi)) \E(\sin(\xi)) \right) ]
\enqs
Note that with this model, the BSDE approach cannot be used  by lack of convexity.

For this example, we take $T=0.4$, $\sigma=1$, 
$\beta(x,\mu)$ $=$  $\bar\mu - x$, and we test the three distributions as given in the Min/max example \ref{sec:MinMax}.

\subsubsection{A two dimensional example}
\label{sec:Systemic2D}
 We consider a multi-dimensional extension of the LQ systemic risk model  of section \ref{sec:Systemic} by supposing that on each dimension, the dynamic satisfies the same equation with independent Brownian motions,  and that the cost functions are 
 the sum over each component of the cost function in the univariate model.  In this case, the value function  is given by  $V(t,x,\mu)$ $=$ $\sum_{i=1}^d V_1(t,x_i,\mu_i)$,  for $t$ $\in$ $[0,T]$, $x$ $=$ $(x_i)_{i\in\llbracket 1,d\rrbracket}$ $\in$ $\R^d$, $\mu_i$ is the $i$-th marginal 
 law of $\mu$ $\in$ $\Pc_2(\R^d)$, and $V_1$ is the value function in the univariate model given by \eqref{eq:vAnal}.

\begin{figure}[H]
\begin{minipage}[t]{0.30\linewidth}
  \centering
 \includegraphics[width=\textwidth]{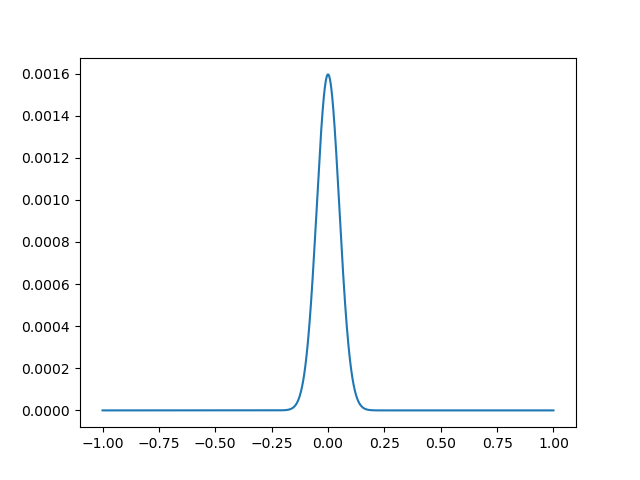}
\caption*{Case 1}
\end{minipage}
\begin{minipage}[t]{0.30\linewidth}
  \centering
 \includegraphics[width=\textwidth]{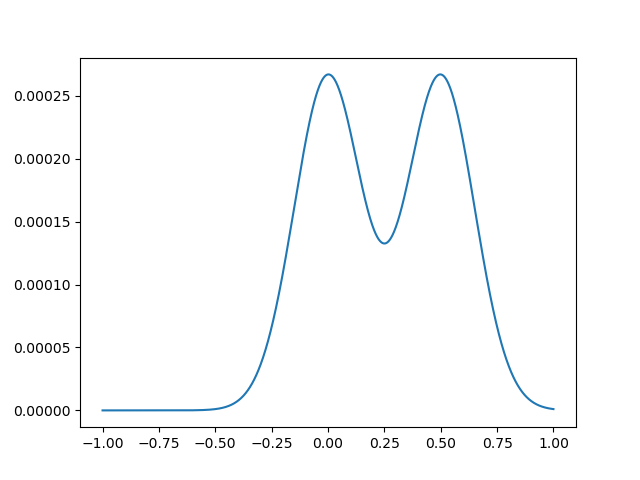}
\caption*{ Case 2}
\end{minipage}
\begin{minipage}[t]{0.30\linewidth}
  \centering
 \includegraphics[width=\textwidth]{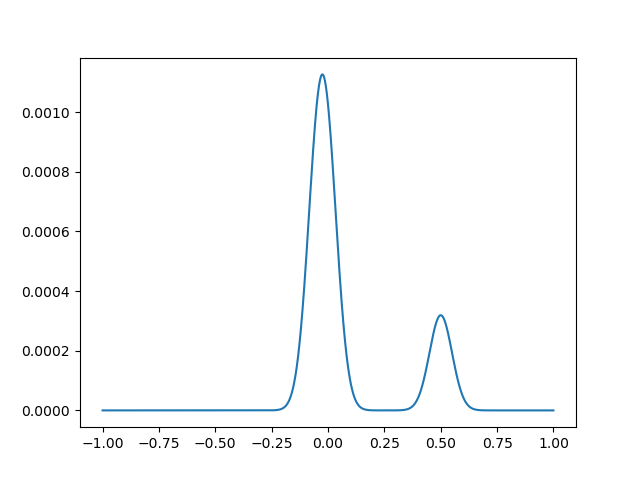}
\caption*{Case 3 }
\end{minipage}
\caption{Marginals of distribution $\mu_0$ tested on the two dimensional systemic model. \label{fig:systemicDist2D} }
 \end{figure}

The parameters of the  dynamic in each dimension  are 
 $\sigma=0.5$, $\kappa =0.6$, $q=0.8$,  $T=0.2$, $c=2$, $\eta=2$. We solve the two dimensional version of the problem \eqref{defV0} by implementing  our various algorithms, and compare  the solution  obtained at $t$ $=$ $0$ with  $v(0,\mu_0)$  for the initial distributions 
$\mu_0$ with the same marginals $\mu_0^1$ $=$ $\mu_0^2$ in the two dimensions 
plotted on Figure \ref{fig:systemicDist2D}:
\begin{itemize}
    \item Case 1 : Gaussian marginals  with $ \bar \mu_0^1 =0$, $std(\mu_0^1)=0.2$,
    \item Case 2 : Mixture of two Gaussian random variables giving the marginal: $X_0^1 =  P (-k_1 + \theta_1 Y) + (1-P) (k_2 + \theta_2 \bar Y)$ with $P$ a  Bernoulli random variable  with parameter $\frac{1}{2}$, $k_1=0$, $k_2= 0.5$, $\theta_1 = \theta_2=  0.15$, $Y, \bar Y \sim  \mathcal{N}(0,1)$,
    \item Case 3 : Mixture of three  Gaussian random variables giving the  marginal: $X_0^1 =   k_{P} + \theta_{P} Y$ where $P$ is a random variable  taking values $(1,2,3)$ with probability $(\frac{2}{5}, \frac{2}{5}, \frac{1}{5})$, $k_1=-0.05$, $k_2=0.$, $k_3=0.5$, $\theta_1 = \theta_2= \theta_3= 0.05$, $Y \sim  \mathcal{N}(0,1)$.
\end{itemize}

\subsection{Results for the systemic risk model 
}

\subsubsection{Dynamic programming-based algorithms}

We report  the results for this model of section \ref{sec:Systemic}  in Tables \ref{tab:globalSys}, \ref{tab:locContSys} and \ref{tab:locSystemic}.  It turns out that the results obtained by Algorithms  \ref{algo:globalcont}   and \ref{alg:fControl} are excellent and very close. We can see that results with $K=100$ or $K=200$ bins for the bins method are very close. Notice that with the bins method, we have to limit the number $K$ of  bins due to memory issues for these two algorithms.
We clearly see the effect of the convergence of the Euler scheme used to discretized the equations on the convergence rate. The Bins method and the Cylinder method provide  very similar results but as the cost of Algorithm \ref{algo:globalcont} is in $O(N_T)$ while the cost of Algorithm \ref{alg:fControl} is in $O(\Frac{N_T(N_T-1)}{2})$, Algorithm \ref{algo:globalcont}  is clearly preferred.
The computational time values presented in Table \ref{tab:locContSys} provide confirmation that Algorithm \ref{alg:fControl} becomes impractical and less usable as the number of time steps increases.

\begin{table}[H]
    \centering
    \begin{tabular}{|c|c|c|c|c|c|c|c|c|c|} \hline
       Method & $K$ &   $\Delta t$ $=$ $\frac{T}{N_T}$ &      \multicolumn{2}{c|}{Case 1} & \multicolumn{2}{c|}{Case 2} & \multicolumn{2}{c|}{Case 3} & Training \\
     & &  &   Calc & Anal & Calc & Anal & Calc & Anal &  time (s) \\ \hline
     Bins & 100 & 0.02  &        0.1670 &  0.1642& 0.1495 &  0.1446  & 0.1497 &   0.1446 & 8160\\
     Bins & 100 & 0.01  &   0.1651&  0.1642& 0.1472 &  0.1446  & 0.1470 &   0.1446 & 16200 \\
     Cylinder & 500 & 0.02 &   0.1684 & 0.1642  & 0.1489 & 0.1446  & 0.1492 & 0.1446  & 8100 \\
     Cylinder & 500 & 0.01 & 0.1665 & 0.1642  & 0.1469 & 0.1446  & 0.1467 & 0.1446 & 15240\\ \hline
    \end{tabular}
    \begin{tabular}{|c|c|c|c|c|c|c|c|c|} \hline
       Method & $K$  &   $\Delta t$ $=$ $\frac{T}{N_T}$ &      \multicolumn{2}{c|}{Case 4} & \multicolumn{2}{c|}{Case 5} & \multicolumn{2}{c|}{Case 6} \\
     & &  &   Calc & Anal & Calc & Anal & Calc & Anal  \\ \hline
      Bins & 100 & 0.02  &       0.1675 & 0.1642 & 0.1824 & 0.1812 & 0.1792  & 0.1772\\
      Bins & 100 & 0.01  &  0.1648 & 0.1642 & 0.1803 & 0.1812 & 0.1766  & 0.1772 \\
      Cylinder & 500 & 0.02 & 0.1684 & 0.1642  & 0.1848 & 0.1812  & 0.1817 & 0.1772 \\
      Cylinder & 500 & 0.01 & 0.1660 & 0.1642  & 0.1835 & 0.1812  & 0.1795 & 0.1772  \\ \hline
    \end{tabular}
    \caption{Global Algorithm \ref{algo:globalcont} for systemic risk with $T=0.2$,  $\Kc = [-1.38, 1.62]$ using $M = 10$, and 60000 gradient iterations. }
    \label{tab:globalSys}
\end{table}

\begin{table}[H]
    \centering
    \begin{tabular}{|c|c|c|c|c|c|c|c|c|c|}  \hline
      Method &  $K$  &    $\Delta t$ $=$  $\frac{T}{N_T}$ &     \multicolumn{2}{c|}{Case 1} & \multicolumn{2}{c|}{Case 2} & \multicolumn{2}{c|}{Case 3} &  Training  \\
     & & & Calc & Anal & Calc & Anal & Calc & Anal  &  time (s) \\  \hline
      Bins & 100 & 0.02  &   0.1692 & 0.1642 &  0.1493 &  0.1446  &  0.1495 & 0.1446 & 20000 \\
      Bins & 100 & 0.01  &   0.1673 & 0.1642  & 0.1478 & 0.1446  & 0.1470 & 0.1446 &  73300 \\
      Bins & 200 & 0.01  &  0.1674 & 0.1642  & 0.1480 & 0.1446  & 0.1477 & 0.1446 &  108800 \\
      Cylinder & 500 & 0.02 & 0.1688 & 0.1642  & 0.1492 & 0.1446  & 0.1490 & 0.1446 & 46600   \\
      Cylinder & 500 & 0.01 & 0.1662 & 0.1642  & 0.1468 & 0.1446  & 0.1471 & 0.1446 &  160200\\ \hline
    \end{tabular}
    \begin{tabular}{|c|c|c|c|c|c|c|c|c|}  \hline
        Method & $K$  & $\Delta t$ $=$  $\frac{T}{N_T}$ &      \multicolumn{2}{c|}{Case 4} & \multicolumn{2}{c|}{Case 5} & \multicolumn{2}{c|}{Case 6} \\
        & &  &  Calc & Anal & Calc & Anal & Calc & Anal  \\  \hline
         Bins & 100 & 0.02  &   0.1691& 0.1642 & 0.1862 & 0.1821 & 0.1822 & 0.1772 \\
         Bins & 100 & 0.01 &  0.1670 & 0.1642  & 0.1836 & 0.1812  & 0.1799 & 0.1772  \\
         Bins & 200 & 0.01  &  0.1675 & 0.1642  & 0.1844 & 0.1812  & 0.1800 & 0.1772 \\
        Cylinder & 500 & 0.02 & 0.1684 & 0.1642  & 0.1862 & 0.1812  & 0.1819 & 0.1772\\
         Cylinder & 500 & 0.01 & 0.1663 & 0.1642  & 0.1836 & 0.1812  & 0.1794 & 0.1772 \\  \hline
    \end{tabular}
    \caption{Policy iteration Algorithm \ref{alg:fControl} for systemic risk with $T=0.2$,  $\Kc = [-1.38, 1.62]$   using $M=10$, and 30000  gradient iterations. }
    \label{tab:locContSys}
\end{table}

The results obtained by the value iteration  Algorithm \ref{alg:LocalfVal} are still  good but less accurate than the results obtained by the two other algorithms. 
The cylinder methods appears to be the best of the two methods. We notice a small degradation of the results as we refine the time step with the bins method. Notice that the memory used by this algorithm is small compared to the two other algorithms and it permits to take a high number $K$ of bins for the bins method (even if it is not necessary on this case).

\begin{table}[H]
 \centering
 \begin{tabular}{|c|c|c|c|c|c|c|c|c|c|}  \hline
 Method & $K$  &  $\Delta t$ $=$  $\frac{T}{N_T}$ &  \multicolumn{2}{c|}{Case 1} & \multicolumn{2}{c|}{Case 2} & \multicolumn{2}{c|}{Case 3} &  Training  \\  \hline
 &  &  & Calc & Anal &  Calc & Anal &  Calc & Anal  & time (s)  \\ 
Bins & 500 &  0.02 &  0.1620 & 0.1642  & 0.1373 & 0.1446  & 0.1698 & 0.1446   &  36000\\
Bins & 500 &  0.01 & 0.1873 & 0.1642  & 0.1673 & 0.1446  & 0.1841 & 0.1446 &  72000\\
Cylinder & 500 & 0.02 & 0.1722 & 0.1642  & 0.1540 & 0.1446  & 0.1554 & 0.1446   & 9300\\
Cylinder & 500 & 0.01 & 0.1704 & 0.1642  & 0.1520 & 0.1446  & 0.1571 & 0.1446 &   18600\\  \hline
 \end{tabular}
 \centering
 \begin{tabular}{|c|c|c|c|c|c|c|c|c|}  \hline
 Method & $K$ &  $\Delta t$ $=$  $\frac{T}{N_T}$ &  \multicolumn{2}{c|}{Case 4 } & \multicolumn{2}{c|}{Case 5} & \multicolumn{2}{c|}{Case 6}  \\
&  &  & Calc & Anal &  Calc & Anal &  Calc & Anal    \\  \hline
Bins & 500 &  0.02 &0.1630 & 0.1642  & 0.1809 & 0.1812  & 0.1755 & 0.1772   \\
Bins & 500 &  0.01 & 0.1880 & 0.1642  & 0.2037 & 0.1812  & 0.1991 & 0.1772   \\
Cylinder & 500 & 0.02 & 0.1722 & 0.1642  & 0.1880 & 0.1812  & 0.1843 & 0.1772   \\
Cylinder & 500 & 0.01 & 0.1704 & 0.1642  & 0.1864 & 0.1812  & 0.1827 & 0.1772   \\  \hline
 \end{tabular}
 \caption{Value iteration Algorithm algorithm  \ref{alg:LocalfVal} for systemic risk with $T=0.2$,  $\Kc = [-1.38, 1.62]$, $M=10$, and 30000 gradient iterations.  \label{tab:locSystemic}}
 \end{table}

 In Table \ref{tab:sensiPD}, we provide sensitivity results for cases 1, 4, and 6 using Algorithms \ref{algo:globalcont} and \ref{alg:LocalfVal} with different methods (bins and cylinder) and parameter settings ($\frac{T}{N_T}=0.02$, $K=100$ for bins, and $K=500$ for cylinder). The results are based on 10 runs, and we report the average value obtained along with the standard deviation. It is observed that the results obtained using different methods and algorithms are generally very similar, except for Algorithm \ref{alg:LocalfVal} with the cylinder network. This particular algorithm shows a higher standard deviation, which is a known characteristic of this approach, as mentioned in \cite{hure2020deep}. Furthermore, all the results seem to converge to the discrete-time solution of the problem, indicating the reliability and accuracy of the algorithms employed.

 \begin{table}[H]
 \centering
 \begin{tabular}{|c|c|c|c|c|c|c|c|c|c|c|}  \hline
  Alg & method &\multicolumn{3}{c|}{Case 1} & \multicolumn{3}{c|}{Case 4} & \multicolumn{3}{c|}{Case 6}  \\  \hline
 & &  Calc & Std & Anal &  Calc & Std & Anal &  Calc & Std & Anal    \\   \hline
\ref{algo:globalcont} &  bins & 0.1691 & 0.0007 & 0.1642  & 0.1692 & 0.0010 & 0.1642  & 0.1825 & 0.0007 & 0.1772   \\  \hline
\ref{algo:globalcont} &  cyl &  0.1687 & 0.0002 & 0.1642  & 0.1686 & 0.0002 & 0.1642 & 0.1816 & 0.0002 & 0.1772   \\  \hline
\ref{alg:fControl} &  bins &  0.1694 & 0.0003 & 0.1642  & 0.1694 & 0.0002 & 0.1642  & 0.1821 & 0.0002 & 0.1772    \\  \hline
\ref{alg:fControl} &  cyl & 0.1687 & 0.0002 & 0.1642  & 0.1687 & 0.0001 & 0.1642  & 0.1815 & 0.0002 & 0.1772    \\  \hline
\ref{alg:LocalfVal} &  bins &  0.1692 & 0.0092 & 0.1642  & 0.1692 & 0.0093 & 0.1642  & 0.1822 & 0.0089 & 0.1772  \\  \hline
\ref{alg:LocalfVal} &  cyl &  0.1807 & 0.0066 & 0.1642  & 0.1807 & 0.0066 & 0.1642  & 0.1931 & 0.0064 & 0.1772    \\  \hline
 \end{tabular}
  \caption{Some sensitivity results using 10 runs. \label{tab:sensiPD}}
 \centering
 \end{table}

 \vspace{3mm}

\subsubsection{Results for Backward SDE-based algorithms}

Results for the systemic example  of section \ref{sec:Systemic} are given in Tables \ref{tab:locBSDESytstemic},  \ref{tab:contBSDESytstemic},  \ref{tab:globalBSDESytstemic}, \ref{tab:sumLocBSDESytstemic} and \ref{tab:globalSumMultiBSDESytstemic}.
All the proposed methods converge very accurately to the solution. As previously seen in the results of the dynamic programming-based algorithms, the number of bins does not  need to be large for the bins network. 
For this test case, the numerical values obtained does not permit to select the best algorithm. 
As Algorithm \ref{alg:BSDElocal2} is by far the most costly,  it should not be the preferred choice.  It is difficult to compare the other algorithms in terms of computing time, but all global algorithms have roughly the same cost in terms of time and  the local deep backward algorithm  \ref{alg:BSDElocal1} is certainly more costly as we have  to achieve an optimization per time step. This drawback due to the number of optimizations  is  reduced by transfer learning, namely  the fact that at each time step the problem is much more smaller to solve as we can initialize the parameters of networks at a given time step by the parameters of  networks of the preceding time step.  
On the other hand, we point out that  all the global algorithms are too far memory consuming to be able to compete with the local deep backward algorithm \ref{alg:BSDElocal1} which seems to be globally the best choice.

\begin{table}[H]
 \centering
 \begin{tabular}{|c|c|c|c|c|c|c|c|c|c|}  \hline
 Method & $K$ &   $\frac{T}{N_T}$ &  \multicolumn{2}{c|}{Case 1} & \multicolumn{2}{c|}{Case 2} & \multicolumn{2}{c|}{Case 3}&  Training   \\
 &  &  & Calc & Anal &  Calc & Anal &  Calc & Anal  & time (s) \\   \hline
Bins & 100 &  0.02 &  0.1690 & 0.1642  & 0.1493 & 0.1446  & 0.1497 & 0.1446 &  7500\\
Bins & 200 &  0.02 &  0.1689 & 0.1642  & 0.1494 & 0.1446  & 0.1494 & 0.1446 & 11700\\
Bins & 100 &  0.01 &  0.1664 & 0.1642  & 0.1474 & 0.1446  & 0.1470 & 0.1446 &  15000 \\
Bins & 200 &  0.01 &  0.1664 & 0.1642  & 0.1472 & 0.1446  & 0.1471 & 0.1446  &  23400\\
Cylinder & 500 & 0.02 &  0.1683 & 0.1642  & 0.1491 & 0.1446  & 0.1492 & 0.1446  & 10500 \\
Cylinder & 500 & 0.01 &  0.1664 & 0.1642  & 0.1472 & 0.1446  & 0.1466 & 0.1446 &  21000\\  \hline
 \end{tabular}
 \centering
 \begin{tabular}{|c|c|c|c|c|c|c|c|c|}
 Method & $K$  &   $\frac{T}{N_T}$ &  \multicolumn{2}{c|}{Case 4 } & \multicolumn{2}{c|}{Case 5} & \multicolumn{2}{c|}{Case 6}  \\
&  &  & Calc & Anal &  Calc & Anal &  Calc & Anal    \\   \hline
Bins & 100 &  0.02 &  0.1690 & 0.1642  & 0.1860 & 0.1812  & 0.1816 & 0.1772 \\
Bins & 200 &  0.02 & 0.1687 & 0.1642  & 0.1853 & 0.1812  & 0.1818 & 0.1772   \\
Bins & 200 &  0.01 &  0.1669 & 0.1642  & 0.1838 & 0.1812  & 0.1801 & 0.1772  \\
Bins & 200 &  0.01 & 0.1666 & 0.1642  & 0.1835 & 0.1812  & 0.1796 & 0.1772   \\
Cylinder & 500 & 0.02 &  0.1683 & 0.1642  & 0.1858 & 0.1812  & 0.1816 & 0.1772  \\
Cylinder & 500 & 0.01 &  0.1665 & 0.1642  & 0.1837 & 0.1812  & 0.1795 & 0.1772  \\  \hline
 \end{tabular}
 \caption{Local deep backward  BSDE Algorithm  \ref{alg:BSDElocal1}, $T = 0.2$, $\Kc= [ -1.38,1.62 ]$, using $M=10$, and $30000$  gradient iterations. \label{tab:locBSDESytstemic} }
 \end{table}
\begin{table}[H]
 \centering
 \begin{tabular}{|c|c|c|c|c|c|c|c|c|c|}  \hline
 Method & $K$ &   $\frac{T}{N_T}$ &  \multicolumn{2}{c|}{Case 1} & \multicolumn{2}{c|}{Case 2} & \multicolumn{2}{c|}{Case 3}  & Training  \\
 &  &  & Calc & Anal &  Calc & Anal &  Calc & Anal  &   time (s) \\   \hline
Bins & 200 &  0.02 &  0.1709 & 0.1642  & 0.1513 & 0.1446  & 0.1516 & 0.1446 &  30000\\
Bins & 200 &  0.01 &   0.1672 & 0.1642  & 0.1479 & 0.1446  & 0.1475 & 0.1446 & 111000 \\
Cylinder & 500 & 0.02 & 0.1688 & 0.1642  & 0.1494 & 0.1446  & 0.1489 & 0.1446  & 20500 \\
Cylinder & 500 & 0.01 &  0.1663 & 0.1642  & 0.1469 & 0.1446  & 0.1472 &  0.1446  & 68800   \\  \hline
 \end{tabular}
 \centering
 \begin{tabular}{|c|c|c|c|c|c|c|c|c|}  \hline
 Method & $K$ &   $\frac{T}{N_T}$ &  \multicolumn{2}{c|}{Case 4 } & \multicolumn{2}{c|}{Case 5} & \multicolumn{2}{c|}{Case 6}  \\
&  &  & Calc & Anal &  Calc & Anal &  Calc & Anal    \\   \hline
Bins & 200 &  0.02 &  0.1711 & 0.1642  & 0.1881 & 0.1812  & 0.1838 & 0.1772 \\
Bins & 200 &  0.01 &   0.1671 & 0.1642  & 0.1845 & 0.1812  & 0.1800 & 0.1772  \\
Cylinder & 500 & 0.02 &  0.1686 & 0.1642  & 0.1855 & 0.1812  & 0.1817 & 0.1772  \\
Cylinder & 500 & 0.01 &  0.1662 & 0.1642  & 0.1834 & 0.1812  & 0.1787 & 0.1772  \\  \hline
 \end{tabular}
 \caption{Deep backward multi-step Algorithm  \ref{alg:BSDElocal2}, $T = 0.2$, $\Kc= [ -1.38,1.62 ]$  using $M=10$, and  $30000$  gradient iterations.  \label{tab:contBSDESytstemic} }
 \end{table}

\begin{table}[H]
 \centering
 \begin{tabular}{|c|c|c|c|c|c|c|c|c|c|}  \hline
 Method & $K$ &   $\frac{T}{N_T}$ &  \multicolumn{2}{c|}{Case 1} & \multicolumn{2}{c|}{Case 2} & \multicolumn{2}{c|}{Case 3} & Training  \\
 &  &  & Calc & Anal &  Calc & Anal &  Calc & Anal  &  time (s)   \\   \hline
Bins & 100 &  0.02 &  0.1691 & 0.1642  & 0.1496 & 0.1446  & 0.1498 & 0.1446  & 4500 \\
Bins & 200 &  0.02 &   0.1691 & 0.1642  & 0.1495 & 0.1446  & 0.1497 & 0.1446 &  6400 \\
Bins & 200 &  0.01 &   0.1663 & 0.1642  & 0.1468 & 0.1446  & 0.1471 & 0.1446 & 12400\\
Cylinder & 500 & 0.02 & 0.1686 & 0.1642  & 0.1491 & 0.1446  & 0.1492 & 0.1446 &  4530\\
Cylinder & 500 & 0.01 &  0.1665 & 0.1642  & 0.1466 & 0.1446  & 0.1466 & 0.1446   & 8400\\  \hline
 \end{tabular}
 \centering
 \begin{tabular}{|c|c|c|c|c|c|c|c|c|}  \hline
 Method & $K$ &   $\frac{T}{N_T}$ &  \multicolumn{2}{c|}{Case 4 } & \multicolumn{2}{c|}{Case 5} & \multicolumn{2}{c|}{Case 6}  \\
&  &  & Calc & Anal &  Calc & Anal &  Calc & Anal    \\   \hline
Bins & 100 &  0.02 &  0.1692 & 0.1642  & 0.1858 & 0.1812  & 0.1815 & 0.1772 \\
Bins & 200 &  0.02 &  0.1694 & 0.1642  & 0.1863 & 0.1812  & 0.1824 & 0.1772 \\
Bins & 200 &  0.01 &  0.1668 & 0.1642  & 0.1838 & 0.1812  & 0.1793 & 0.1772 \\
Cylinder & 500 & 0.02 &  0.1686 & 0.1642  & 0.1857 & 0.1812  & 0.1816 & 0.1772  \\
Cylinder & 500 & 0.01 &  0.1667 & 0.1642  & 0.1836 & 0.1812  & 0.1795 & 0.1772  \\  \hline
 \end{tabular}
 \caption{Global  deep MKV  BSDE Algorithm \ref{algo:BSDEglobal},  $T = 0.2$,  $\Kc$ $=$ $[-1.38,1.6 ]$  using $M=10$, and  $30000$  gradient iterations.  \label{tab:globalBSDESytstemic} }
 \end{table}

\begin{table}[H]
 \centering
 \begin{tabular}{|c|c|c|c|c|c|c|c|c|c|}  \hline
 Method & $K$ &   $\frac{T}{N_T}$ &  \multicolumn{2}{c|}{Case 1} & \multicolumn{2}{c|}{Case 2} & \multicolumn{2}{c|}{Case 3}  &  Training  \\
 &  &  & Calc & Anal &  Calc & Anal &  Calc & Anal   &  time (s)  \\   \hline
Bins & 200 &  0.02 &  0.1753 & 0.1642  & 0.1545 & 0.1446  & 0.1706 & 0.1446 &  6400 \\
Bins & 200 &  0.01 &  0.1670 & 0.1642  & 0.1483 & 0.1446  & 0.1597 & 0.1446  & 12300 \\
Cylinder & 500 & 0.02 & 0.1684 & 0.1642  & 0.1496 & 0.1446  & 0.1491 & 0.1446  & 4200  \\
Cylinder & 500 & 0.01 &  0.1667 & 0.1642  & 0.1469 & 0.1446  & 0.1468 & 0.1446 &  8100  \\  \hline
 \end{tabular}
 \centering
 \begin{tabular}{|c|c|c|c|c|c|c|c|c|}  \hline
 Method & $K$ &   $\frac{T}{N_T}$ &  \multicolumn{2}{c|}{Case 4 } & \multicolumn{2}{c|}{Case 5} & \multicolumn{2}{c|}{Case 6}  \\
&  &  & Calc & Anal &  Calc & Anal &  Calc & Anal    \\   \hline
Bins & 200 &  0.02 & 0.1758 & 0.1642  & 0.1931 & 0.1812  & 0.1887 & 0.1772   \\
Bins & 200 &  0.01 & 0.1661 & 0.1642  & 0.1841 & 0.1812  & 0.1797 & 0.1772  \\
Cylinder & 500 & 0.02 &  0.1687 & 0.1642  & 0.1856 & 0.1812  & 0.1816 & 0.1772  \\
Cylinder & 500 & 0.01 &   0.1664 & 0.1642  & 0.1836 & 0.1812  & 0.1793 & 0.1772 \\  \hline
 \end{tabular}
 \caption{Global/local deep MKV   BSDE Algorithm  \ref{algo:BSDEglobal2},  $T = 0.2$,  $\Kc$ $=$ $[-1.38,1.62 ]$ using $M=10$, $30000$  gradient iterations.  \label{tab:sumLocBSDESytstemic} }
 \end{table}

 \begin{table}[H]
 \centering
 \begin{tabular}{|c|c|c|c|c|c|c|c|c|c|}  \hline
 Method & $K$ &   $\frac{T}{N_T}$ &  \multicolumn{2}{c|}{Case 1} & \multicolumn{2}{c|}{Case 2} & \multicolumn{2}{c|}{Case 3} & Training  \\ 
 &  &  & Calc & Anal &  Calc & Anal &  Calc & Anal  & time (s)  \\  \hline
Bins & 200 &  0.02 &   0.1689 & 0.1642  & 0.1507 & 0.1446  & 0.1528 & 0.1446 &   6300 \\
Bins & 200 &  0.01 &  0.1664 & 0.1642  & 0.1470 & 0.1446  & 0.1469 & 0.1446 &   12400 \\
Cylinder & 500 & 0.02 &   0.1685 & 0.1642  & 0.1489 & 0.1446  & 0.1494 & 0.1446 & 4300 \\
Cylinder & 500 & 0.01 &   0.1658 & 0.1642  & 0.1470 & 0.1446  & 0.1468 & 0.1446 &  83200 \\  \hline
 \end{tabular}
 \centering
 \begin{tabular}{|c|c|c|c|c|c|c|c|c|}  \hline
 Method & $K$ &   $\frac{T}{N_T}$ &  \multicolumn{2}{c|}{Case 4 } & \multicolumn{2}{c|}{Case 5} & \multicolumn{2}{c|}{Case 6}  \\
&  &  & Calc & Anal &  Calc & Anal &  Calc & Anal    \\   \hline
Bins & 200 &  0.02 & 0.1692 & 0.1642  & 0.1868 & 0.1812  & 0.1821 & 0.1772   \\
Bins & 200 &  0.01 &  0.1666 & 0.1642  & 0.1829 & 0.1812  & 0.1796 & 0.1772 \\
Cylinder & 500 & 0.02 &  0.1687 & 0.1642  & 0.1855 & 0.1812  & 0.1817 & 0.1772  \\
Cylinder & 500 & 0.01 & 0.1661 & 0.1642  & 0.1834 & 0.1812  & 0.1795 & 0.1772   \\  \hline
 \end{tabular}
 \caption{Global deep multi-step MKV BSDE Algorithm  \ref{algo:BSDEglobalmulti}, $T = 0.2$, $\Kc= [-1.38,1.62]$  using $M=10$, $30000$  gradient iterations.  \label{tab:globalSumMultiBSDESytstemic} }
 \end{table}

In Table \ref{tab:sensiBSDEPD}, sensitivity results are presented based on ten runs using the same hyperparameters as in the dynamic programming approach. Similar to the previous table, all the algorithms demonstrate consistent results with low standard deviations, except for Algorithm \ref{algo:BSDEglobalmulti} with bins, which exhibits a higher standard deviation. Additionally, it is observed that all the algorithms converge to the same value as in the dynamic programming approach, further confirming their reliability and accuracy.

 \begin{table}[H]
 \centering
 \begin{tabular}{|c|c|c|c|c|c|c|c|c|c|c|}  \hline
  Alg & method &\multicolumn{3}{c|}{Case 1} & \multicolumn{3}{c|}{Case 4} & \multicolumn{3}{c|}{Case 6}  \\  \hline
 & &  Calc & Std & Anal &  Calc & Std & Anal &  Calc & Std & Anal    \\   \hline
\ref{alg:BSDElocal1} &  bins &  0.1691 & 0.0002 & 0.1642  & 0.1691 & 0.0003 & 0.1642  & 0.1821 & 0.0003 & 0.1772  \\  \hline
\ref{alg:BSDElocal1} &  cyl &  0.1688 & 0.0001 & 0.1642  & 0.1687 & 0.0002 & 0.1642  & 0.1817 & 0.0002 & 0.1772   \\  \hline
\ref{alg:BSDElocal2} &  bins &   0.1691 & 0.0002 &  0.1642   & 0.1691 & 0.0002 &  0.1642  &  0.1820 & 0.0002 & 0.1772   \\  \hline
\ref{alg:BSDElocal2} &  cyl &   0.1686 & 0.0003 & 0.1642  & 0.1689 & 0.0003 & 0.1642  & 0.1816 & 0.0001 & 0.1772 \\  \hline
\ref{algo:BSDEglobal} &  bins &  0.1693 & 0.0005 & 0.1642  & 0.1693 & 0.0004 & 0.1642  & 0.1821 & 0.0004 & 0.1772   \\  \hline
\ref{algo:BSDEglobal} &  cyl &  0.1687 & 0.0002 & 0.1642  & 0.1689 & 0.0003 & 0.1642  & 0.1817 & 0.0002 & 0.1772  \\  \hline
\ref{algo:BSDEglobal2} &  bins & 0.1709 & 0.0013 & 0.1642  & 0.1704 & 0.0013 & 0.1642  & 0.1830 & 0.0010 & 0.1772   \\  \hline
\ref{algo:BSDEglobal2} &  cyl &  0.1686 & 0.0002 & 0.1642  & 0.1686 & 0.0003 & 0.1642  & 0.1817 & 0.0002 & 0.1772  \\  \hline
\ref{algo:BSDEglobalmulti} &  bins &  0.1691 & 0.0004 & 0.1642  & 0.1691 & 0.0003 & 0.1642  & 0.1820 & 0.0004 & 0.1772  \\  \hline
\ref{algo:BSDEglobalmulti} &  cyl &  0.1687 & 0.0002 & 0.1642  & 0.1687 & 0.0002 & 0.1642  & 0.1815 & 0.0002 & 0.1772  \\  \hline
 \end{tabular}
  \caption{Some sensitivity results using 10 runs. \label{tab:sensiBSDEPD}}
 \centering
 \end{table}

 \subsection{Results for the min/max MKV model}
 
 \subsubsection{Dynamic programming-based algorithms}

 Results for $T=0.2$ are reported in table \ref{tab:MinDisGlobalT02}, \ref{tab:MinDisContT02}, \ref{tab:MinDisLocT02}: they  are very good  for all algorithms and network used.
 Results for $T=0.5$ are  reported in table \ref{tab:MinDisGlobalT05}, \ref{tab:MinDisContT05}, \ref{tab:MinDisLocT05}, and   also give excellent results. 
 Notice that with  Algorithm \ref{alg:fControl},  it is impossible to solve the problem with $T=0.5$ using $N_T=50$ due to memory issues and the time needed limited to 3 days.
 As we increase the number of time steps for Algorithm \ref{alg:LocalfVal}, we observe for the bins methods, as in the previous test case,  a small degradation of the results due to an accumulation of regression error, and therefore Algorithm \ref{algo:globalcont} should be preferred.  
 
 \begin{table}[H]
    \centering
    \begin{tabular}{|c|c|c|c|c|c|c|c|c|}  \hline
        &   & & \multicolumn{2}{c|}{Case 1} & \multicolumn{2}{c|}{Case 2} &  \multicolumn{2}{c|}{Case 3}  \\
      Method  & $K$ &  $\frac{T}{N_T}$  &  Calc & Ref & Calc & Ref & Calc & Ref \\  \hline
        Bins & 100 & 0.02 & 0.481 & 0.483 & 0.502 & 0.494 &  0.489 &   0.491 \\ 
        Bins & 100 & 0.01 & 0.481 & 0.483 & 0.503 & 0.494 &  0.489 &   0.491  \\
        Bins & 200 & 0.01 & 0.484 & 0.483  & 0.498 & 0.494  & 0.491 & 0.491 \\ 
        Cylinder & 500 & 0.02 &  0.484 & 0.483  & 0.493 & 0.494  & 0.491 & 0.491   \\
       Cylinder & 500 & 0.01 &   0.484 & 0.483  & 0.494 & 0.494  & 0.491 & 0.491 \\  \hline
    \end{tabular}
    \caption{Global Algorithm \ref{algo:globalcont} with   $T=0.2$,  $ \Kc= [0.21,2.72]$. }
    \label{tab:MinDisGlobalT02}
\end{table}

\begin{table}[H]
    \centering
    \begin{tabular}{|c|c|c|c|c|c|c|c|c|}  \hline
        &   & & \multicolumn{2}{c|}{Case 1} & \multicolumn{2}{c|}{Case 2} &  \multicolumn{2}{c|}{Case 3}  \\
      Method  & $K$ &  $\frac{T}{N_T}$  &  Calc & Ref & Calc & Ref & Calc & Ref \\  \hline
        Bins & 100 & 0.02 &   0.830  &  0.818 & 1.100 &1.082  & 0.848  &   0.836 \\
        Bins & 100 & 0.01 &   0.833  &  0.818 & 1.104 &1.082  & 0.850  &   0.836\\
        Bins & 200 & 0.01 &  0.831  &  0.818 & 1.092 &1.082  & 0.848  &   0.836   \\ 
        Cylinder & 500 & 0.02 & 0.814 & 0.818  & 1.080 & 1.082  & 0.831 & 0.836   \\
        Cylinder & 500 & 0.01&  0.819 & 0.818  & 1.085 & 1.082  & 0.837 & 0.836   \\  \hline
    \end{tabular}
    \caption{Global Algorithm \ref{algo:globalcont} with  $T=0.5$, $ \Kc = [-0.4 , 3.21]$. }
    \label{tab:MinDisGlobalT05}
\end{table}

\begin{table}[H]
    \centering
    \begin{tabular}{|c|c|c|c|c|c|c|c|c|}  \hline
        &   & & \multicolumn{2}{c|}{Case 1} & \multicolumn{2}{c|}{Case 2} &  \multicolumn{2}{c|}{Case 3}  \\
      Method  & $K$  &  $\frac{T}{N_T}$  &  Calc & Ref & Calc & Ref & Calc & Ref \\  \hline
        Bins & 100 & 0.02 & 0.480 & 0.483  & 0.502 & 0.494  & 0.489 & 0.491 \\ 
        Bins & 200 & 0.02 & 0.482 & 0.483  & 0.496 & 0.494  & 0.491 & 0.491 \\ 
        Cylinder & 500 & 0.02 &   0.484 & 0.483  & 0.493 & 0.494  & 0.491 & 0.491  \\  \hline
    \end{tabular}
    \caption{Policy iteration Algorithm \ref{alg:fControl}  with   $T=0.2$,  $ \Kc= [0.21,2.72]$.}
    \label{tab:MinDisContT02}
\end{table}

\begin{table}[H]
    \centering
    \begin{tabular}{|c|c|c|c|c|c|c|c|c|}  \hline
        &   & & \multicolumn{2}{c|}{Case 1} & \multicolumn{2}{c|}{Case 2} &  \multicolumn{2}{c|}{Case 3}  \\
      Method  & $K$ &  $\frac{T}{N_T}$  &  Calc & Ref & Calc & Ref & Calc & Ref \\  \hline
        Bins & 100 & 0.02 &   0.819 & 0.818  & 1.088 & 1.082  & 0.836 & 0.836  \\
         Bins & 200 & 0.02 & 0.818 & 0.818  & 1.090 & 1.082  & 0.836 & 0.836  \\
         Cylinder & 500 & 0.02 &  0.814 & 0.818  & 1.081 & 1.082  & 0.831 & 0.836  \\  \hline
    \end{tabular}
    \caption{Policy iteration Algorithm \ref{alg:fControl} with  $T=0.5$, $ \Kc = [-0.4 , 3.21]$.}
    \label{tab:MinDisContT05}
\end{table}

\begin{table}[H]
    \centering
    \begin{tabular}{|c|c|c|c|c|c|c|c|c|}  \hline
        &   & & \multicolumn{2}{c|}{Case 1} & \multicolumn{2}{c|}{Case 2} &  \multicolumn{2}{c|}{Case 3}  \\
      Method  & $K$ &  $\frac{T}{N_T}$  &  Calc & Ref & Calc & Ref & Calc & Ref \\  \hline
        Bins & 100 & 0.02 &   0.494 & 0.483 & 0.512 & 0.494 &  0.502 &   0.491 \\ 
        Bins & 200 & 0.02 &  0.490 & 0.483  & 0.493 & 0.494  & 0.495 & 0.491   \\
        Cylinder & 500 & 0.02 &   0.486 & 0.483  & 0.493 & 0.494  & 0.491 & 0.491 \\ \hline
    \end{tabular}
    \caption{Value iteration Algorithm \ref{alg:LocalfVal}  with   $T=0.2$,  $ \Kc= [0.21,2.72]$. }
    \label{tab:MinDisLocT02}
\end{table}

\begin{table}[H]
    \centering
    \begin{tabular}{|c|c|c|c|c|c|c|c|c|}  \hline
        &   & & \multicolumn{2}{c|}{Case 1} & \multicolumn{2}{c|}{Case 2} &  \multicolumn{2}{c|}{Case 3}  \\
      Method  & $K$ &  $\frac{T}{N_T}$  &  Calc & Ref & Calc & Ref & Calc & Ref \\  \hline
        Bins & 100 & 0.02 &   0.800  &  0.818 & 1.084 &1.082  &  0.817 &   0.836 \\
         Bins & 200 & 0.02 & 0.810 & 0.818  & 1.079 & 1.082  & 0.828 & 0.836 \\
                 Bins & 200 & 0.01 & 0.835 & 0.818  & 1.114 & 1.082  & 0.853 & 0.836  \\
         Cylinder & 500 & 0.02 & 0.811 & 0.818  & 1.088 & 1.082  & 0.829 & 0.836  \\
        Cylinder & 500 & 0.01 &  0.810 & 0.818  & 1.078 & 1.082  & 0.827 & 0.836  \\  \hline
    \end{tabular}
    \caption{Value iteration Algorithm \ref{alg:LocalfVal} with  $T=0.5$, $ \Kc = [-0.4 , 3.21]$.}
    \label{tab:MinDisLocT05}
\end{table}
It is important to consider that as the maturity increases, the size of $\Kc$ needs to be adjusted accordingly to ensure that the particles primarily remain within $\Kc$. This adjustment is necessary to accommodate the potential expansion of the distribution's support as the maturity lengthens.

\subsubsection{Results for Backward SDE-based algorithms}

Results for this example  of Section \ref{sec:MinMax}  are reported  in Tables \ref{tab:LocBSDEMinDist}, \ref{tab:controlBSDEMinDist}, \ref{tab:globalBSDEMinDist}, \ref{tab:globalSumLocBSDEMinDist} and \ref{tab:globalSumStepBSDEMinDist}. 
All algorithms seem to converge to the good solution except the global deep MKV BSDE Algorithm \ref{algo:BSDEglobal} that always converges on our tests (repeated many times)  to a slightly different solution while using the cylinder network. Notice that, by using the bins network, we avoid the problem on this test case. 
Again it is not feasible to refine the time step when  implementing  the deep backward multi-step Algorithm \ref{alg:BSDElocal2} due to the computational time taken by the algorithm.  The local deep backward  Algorithm  \ref{alg:BSDElocal1}  seems to be the best as the results obtained in Table \ref{tab:LocBSDEMinDist} 
are very good and the memory needed rather small. Either bins or cylinder networks can be used. 
 
\begin{table}[H]
 \centering
 \begin{tabular}{|c|c|c|c|c|c|c|c|c|}  \hline
 Method & $K$ &   $\frac{T}{N_T}$ &  \multicolumn{2}{c|}{Case 1} & \multicolumn{2}{c|}{Case 2} & \multicolumn{2}{c|}{Case 3}  \\
 &  &  & Calc & Ref &  Calc & Ref &  Calc & Ref    \\   \hline
Bins & 100 &  0.02 & 0.8355 & 0.8180  & 1.1074 & 1.0820  & 0.8537 & 0.8360   \\
Bins & 200 &  0.02 &  0.8278 & 0.8180  & 1.0962 & 1.0820  & 0.8462 & 0.8360     \\
Bins & 200 &  0.01 &   0.8343 & 0.8180  & 1.0998 & 1.0820  & 0.8513 & 0.8360   \\
Cylinder & 500 & 0.02 &  0.8249 & 0.8180  & 1.0896 & 1.0820  & 0.8427 & 0.8360  \\
Cylinder & 500 & 0.01 &  0.8312 & 0.8180  & 1.0946 & 1.0820  & 0.8487 & 0.8360  \\  \hline
 \end{tabular}
 \caption{Deep backward  Algorithm  \ref{alg:BSDElocal1}, $T = 0.5$, $\Kc= [-0.40,3.21 ]$.\label{tab:LocBSDEMinDist} }
 \end{table}
\begin{table}[H]
 \centering
 \begin{tabular}{|c|c|c|c|c|c|c|c|c|}  \hline
 Method & $K$ &   $\frac{T}{N_T}$ &  \multicolumn{2}{c|}{Case 1} & \multicolumn{2}{c|}{Case 2} & \multicolumn{2}{c|}{Case 3}  \\
 &  &  & Calc & Ref &  Calc & Ref &  Calc & Ref    \\   \hline
Bins & 200 &  0.02 & 0.8277 & 0.8180  & 1.0966 & 1.0820  & 0.8453 & 0.8360      \\
Cylinder & 500 & 0.02 & 0.8259 & 0.8180  & 1.0904 & 1.0820  & 0.8427 & 0.8360 \\  \hline
 \end{tabular}
 \caption{Deep backward multi-step  Algorithm  \ref{alg:BSDElocal2}, $T = 0.5$, $\Kc= [-0.40,3.21]$. \label{tab:controlBSDEMinDist} }
 \end{table}
 \begin{table}[H]
 \centering
 \begin{tabular}{|c|c|c|c|c|c|c|c|c|}  \hline
 Method & $K$ &   $\frac{T}{N_T}$ &  \multicolumn{2}{c|}{Case 1} & \multicolumn{2}{c|}{Case 2} & \multicolumn{2}{c|}{Case 3}  \\
 &  &  & Calc & Ref &  Calc & Ref &  Calc & Ref    \\   \hline
Bins & 200 &  0.02 &  0.8299 & 0.8180  & 1.0977 & 1.0820  & 0.8485 & 0.8360    \\
Bins & 100 &  0.01 &  0.8447 & 0.8180  & 1.1111 & 1.0820  & 0.8625 & 0.8360    \\
Bins & 200 &  0.01 &  0.8369 & 0.8180  & 1.1018 & 1.0820  & 0.8566 & 0.8360     \\
Cylinder & 500 & 0.02 & 0.7801 & 0.8180  & 1.0493 & 1.0820  & 0.7968 & 0.8360   \\
Cylinder & 500 & 0.01 & 0.7597 & 0.8180  & 1.0325 & 1.0820  & 0.7767 & 0.8360  \\  \hline
 \end{tabular}
 \caption{Global deep MKV BSDE Algorithm  \ref{algo:BSDEglobal}, $T = 0.5$, $\Kc=[-0.40,3.21]$. \label{tab:globalBSDEMinDist} }
 \end{table}
 \begin{table}[H]
 \centering
 \begin{tabular}{|c|c|c|c|c|c|c|c|c|}  \hline
 Method & $K$ &   $\frac{T}{N_T}$ &  \multicolumn{2}{c|}{Case 1} & \multicolumn{2}{c|}{Case 2} & \multicolumn{2}{c|}{Case 3}  \\
 &  &  & Calc & Ref &  Calc & Ref &  Calc & Ref    \\   \hline
Bins & 200 &  0.02 &  0.8528 & 0.8180  & 1.1003 & 1.0820  & 0.8692 & 0.8360   \\
Bins & 100 &  0.01 &  0.9146 & 0.8180  & 1.1120 & 1.0820  & 0.9219 & 0.8360    \\
Bins & 200 &  0.01 &  0.8406 & 0.8180  & 1.1001 & 1.0820  & 0.8560 & 0.8360   \\
Cylinder & 500 & 0.02 &  0.8305 & 0.8180  & 1.0952 & 1.0820  & 0.8466 & 0.8360   \\
Cylinder & 500 & 0.01 &  0.8666 & 0.8180  & 1.1104 & 1.0820  & 0.8817 & 0.8360  \\  \hline
 \end{tabular}
 \caption{Global/local deep MKV  BSDE Algorithm  \ref{algo:BSDEglobal2}, $T = 0.5$, $\Kc= [  -0.40  , 3.21 ]$. \label{tab:globalSumLocBSDEMinDist} }
 \end{table}
 \begin{table}[H]
 \centering
 \begin{tabular}{|c|c|c|c|c|c|c|c|c|}  \hline
 Method & $K$ &   $\frac{T}{N_T}$ &  \multicolumn{2}{c|}{Case 1} & \multicolumn{2}{c|}{Case 2} & \multicolumn{2}{c|}{Case 3}  \\
 &  &  & Calc & Ref &  Calc & Ref &  Calc & Ref    \\   \hline
Bins & 200 & 0.02 & 0.8380 & 0.8180  & 1.1004 & 1.0820  & 0.8497 & 0.8360  \\
Bins & 200 &  0.01 &    0.8353 & 0.8180  & 1.1002 & 1.0820  & 0.8520 & 0.8360 \\
Cylinder & 500 & 0.02 &  0.8265 & 0.8180  & 1.0902 & 1.0820  & 0.8434 & 0.8360 \\
Cylinder & 500 & 0.01 &   0.8319 & 0.8180  & 1.0951 & 1.0820  & 0.8487 & 0.8360 \\  \hline
 \end{tabular}
\caption{Global deep multi-step MKV BSDE Algorithm  \ref{algo:BSDEglobalmulti},  $T = 0.5$, $\Kc= [  -0.40  , 3.21 ]$.  \label{tab:globalSumStepBSDEMinDist} }
 \end{table}

\subsection{Result on the mean variance problem using the dynamic programming approach.}

We do not report results from  Algorithm \ref{alg:LocalfVal}: indeed, they  diverge for all discretizations tested. Results for the two other algorithms are given in Tables \ref{tab:globalPortT02} and \ref{tab:contPortT02} for $T$ $=$ $0.2$, and 
in Tables  \ref{tab:globalPortT05} and \ref{tab:contPortT05} for $T$ $=$ $0.5$.  Notice that the number of bins taken for the bins network has to be high   to get an accurate solution.

\begin{table}[H]
    \centering
    \begin{tabular}{|c|c|c|c|c|c|c|c|}  \hline
       Method  & $K$ &   \multicolumn{2}{c|}{Case 1} & \multicolumn{2}{c|}{Case 2} &  \multicolumn{2}{c|}{Case 3} \\
      &&  Calc & Anal & Calc & Anal & Calc & Anal   \\  \hline
      Bins & 100 &      -0.0954 &  -0.0865 & -0.1147  &  -0.1059  &  -0.3139 &   -0.3050  \\
       Bins & 200 &    -0.0907 &  -0.0865 & -0.1104  &  -0.1059  &  -0.3094 &   -0.3050  \\
       Bins & 400 &   -0.0882 &  -0.0865 & -0.1081  &  -0.1059  &  -0.3071 &   -0.3050  \\
       Cylinder & 500 &-0.0884 & -0.0865  & -0.1078 & -0.1060  & -0.3070 & -0.3051  \\  \hline
    \end{tabular}
        \begin{tabular}{|c|c|c|c|c|c|c|c|}  \hline
         Method  & $K$ &   \multicolumn{2}{c|}{Case 4} & \multicolumn{2}{c|}{Case 5} & \multicolumn{2}{c|}{Case 6} \\
         &&  Calc & Anal & Calc & Anal & Calc & Anal  \\  \hline
           Bins & 100 &  -0.0952 & -0.0865 & -0.0547& -0.0464 & -0.1769 & -0.1683 \\
       Bins & 200 &    -0.0908 & -0.0865 & -0.0510 & -0.0464 & -0.1724 & -0.1683 \\
       Bins & 400 &  -0.0894 & -0.0865 & -0.0487 & -0.0464 &   -0.1703 & -0.1683 \\ 
       Cylinder & 500 & -0.0883 & -0.0865  & -0.0485 & -0.0464  & -0.1703 & -0.1683  \\  \hline
    \end{tabular}
    \caption{Global Algorithm \ref{algo:globalcont}, $T=0.2$, $ \Kc = [-0.85, 0.9]$, $\frac{T}{N_T} =0.02$.}
    \label{tab:globalPortT02}
\end{table}

\begin{table}[H]
    \centering
    \begin{tabular}{|c|c|c|c|c|c|c|c|}  \hline
       Method  & $K$ &   \multicolumn{2}{c|}{Case 1} & \multicolumn{2}{c|}{Case 2} &  \multicolumn{2}{c|}{Case 3} \\
      &&   Calc & Anal & Calc & Anal & Calc & Anal   \\  \hline
       Bins & 200 &  -0.1018 & -0.0965  & -0.1214 & -0.1156  & -0.3200 & -0.3147  \\
       Bins & 400 &  -0.0976 & -0.0965  & -0.1163 & -0.1156  & -0.3149 & -0.3147 \\
       Cylinder & 500 & -0.0987 & -0.0965  & -0.1179 & -0.1156  & -0.3172 & -0.3147   \\  \hline
    \end{tabular}
        \begin{tabular}{|c|c|c|c|c|c|c|c|}  \hline
         Method  & $K$ & \multicolumn{2}{c|}{Case 4} & \multicolumn{2}{c|}{Case 5} & \multicolumn{2}{c|}{Case 6} \\
         &&  Calc & Anal & Calc & Anal & Calc & Anal  \\  \hline
       Bins & 200 &   -0.1022 & -0.0965  & -0.0613 & -0.0562  & -0.1842 & -0.1786   \\
       Bins & 400 &  -0.0969 & -0.0965  & -0.0562 & -0.0562  & -0.1788 & -0.1786   \\ 
       Cylinder & 500 &  -0.0985 & -0.0965  & -0.0583 & -0.0562  & -0.1804 & -0.1786 \\ \hline
    \end{tabular}
    \caption{Global Algorithm \ref{algo:globalcont}, $T=0.5$, $ \Kc = [-0.85, 0.9]$, $\frac{T}{N_T} =0.02$.}
    \label{tab:globalPortT05}
\end{table}

\begin{table}[H]
 \centering
 \begin{tabular}{|c|c|c|c|c|c|c|c|}  \hline
 Method & $K$ &  \multicolumn{2}{c|}{Case 1} & \multicolumn{2}{c|}{Case 2} & \multicolumn{2}{c|}{Case 3}  \\
 &  &   Calc & Anal &  Calc & Anal &  Calc & Anal    \\   \hline
Bins & 100 &   -0.0959 & -0.0865  & -0.1143 & -0.1060  & -0.3138 & -0.3051   \\
 Bins & 200 &   -0.0906 &  -0.0865 & -0.1102 &  -0.1059  &  -0.3094 &   -0.3050  \\
 Bins & 400 &    -0.0884 &  -0.0865 &  -0.1083 &  -0.1059  & -0.3072  &   -0.3050  \\
 Cylinder & 500 & -0.0884 & -0.0865  & -0.1078 & -0.1060  & -0.3070 & -0.3051    \\  \hline
 \end{tabular}
 \centering
 \begin{tabular}{|c|c|c|c|c|c|c|c|}  \hline
 Method  & $K$ &   \multicolumn{2}{c|}{Case 4 } & \multicolumn{2}{c|}{Case 5} & \multicolumn{2}{c|}{Case 6}  \\
 &  &   Calc & Anal &  Calc & Anal &  Calc & Anal    \\   \hline
Bins & 100 &  -0.0954 & -0.0865  & -0.0553 & -0.0464  & -0.1766 & -0.1683   \\
 Bins & 200 &  -0.0908 & -0.0865 & -0.0505 & -0.0464 &   -0.1723 & -0.1683 \\
 Bins & 400 &      -0.0887 & -0.0865 & -0.0482 & -0.0464 &   -0.1704 & -0.1683 \\ 
 Cylinder & 500 &   -0.0883 & -0.0865  & -0.0485 & -0.0464  & -0.1703 & -0.1683 \\  \hline
 \end{tabular}
 \caption{Policy iteration Algorithm   \ref{alg:fControl}  , $T=0.2$, $ \Kc = [-0.85, 0.9]$, $\frac{T}{N_T} =0.02$.}
\label{tab:contPortT02}
 \end{table}

\begin{table}[H]
 \centering
 \begin{tabular}{|c|c|c|c|c|c|c|c|}  \hline
 Method & $K$ &   \multicolumn{2}{c|}{Case 1} & \multicolumn{2}{c|}{Case 2} & \multicolumn{2}{c|}{Case 3}  \\
 &  &  Calc & Anal &  Calc & Anal &  Calc & Anal    \\   \hline
 Bins & 400 & -0.0978 & -0.0965  & -0.1171 & -0.1156  & -0.3140 & -0.3147   \\  
 Cylinder & 500 &  -0.0986 & -0.0965  & -0.1175 & -0.1156  & -0.3164 & -0.3147   \\  \hline
 \end{tabular}
 \centering
 \begin{tabular}{|c|c|c|c|c|c|c|c|}  \hline
 Method  & $K$  &   \multicolumn{2}{c|}{Case 4 } & \multicolumn{2}{c|}{Case 5} & \multicolumn{2}{c|}{Case 6}  \\
 &  &  Calc & Anal &  Calc & Anal &  Calc & Anal    \\   \hline
 Bins & 400 &   -0.0985 & -0.0965  & -0.0579 & -0.0562  & -0.1789 & -0.1786  \\ 
 Cylinder & 500 & -0.0986 & -0.0965  & -0.0583 & -0.0562  & -0.1807 & -0.1786    \\  \hline
 \end{tabular}
 \caption{Policy iteration Algorithm   \ref{alg:fControl}  , $T=0.5$, $ \Kc = [-0.85, 0.9]$, $\frac{T}{N_T} =0.02$.
}
\label{tab:contPortT05}
 \end{table}

 Again, in term of accuracy, Algorithms \ref{algo:globalcont} and  \ref{alg:fControl} give similar accurate results and the memory taken by both algorithms is close. However Algorithm \ref{algo:globalcont} has to be preferred as the computation time is far lower when we are interested by computing the solution only at time $t$ $=$ $0$. 


\subsection{Results for the  non LQ MKV  model using dynamic programming}

In Table \ref{tab:NLQ1}, results from different algorithms are provided for one run with a maturity $T=0.4$. The settings used include $N_T=20$, $30000$ gradient iterations, $M=10$, and $\Kc= [ -0.81, 2.81 ]$.
Based on the provided results, it is observed that Algorithm \ref{alg:LocalfVal} yields inaccurate results. However, Algorithms \ref{algo:globalcont} and \ref{alg:fControl} demonstrate a  high accuracy in capturing the desired outcome.

\begin{table}[H]
 \centering
 \begin{tabular}{|c|c|c|c|c|c|c|c|}  \hline
  Alg & method & \multicolumn{2}{c|}{Case 1 }& \multicolumn{2}{c|}{Case 2} & \multicolumn{2}{c|}{Case 3}  \\  \hline
 & &  Calc & Anal &  Calc & Anal &  Calc & Anal    \\   \hline
\ref{algo:globalcont} &  bins &  1.1840 & 1.1735  & 0.9389 & 0.9200  & 1.1695 & 1.1585  \\  \hline
\ref{algo:globalcont} &  cyl &  1.1815 & 1.1736  & 0.9233 & 0.9197  & 1.1666 & 1.1585  \\  \hline
\ref{alg:fControl} &  bins & 1.1896 & 1.1735  & 0.9479 & 0.9197  & 1.1757 & 1.1584   \\  \hline
\ref{alg:fControl} &  cyl & 1.1791 & 1.1735  & 0.9239 & 0.9196  & 1.1640 & 1.1585    \\  \hline
\ref{alg:LocalfVal} &  bins &  0.9356 & 1.1736  & 0.8282 & 0.9198  & 0.9293 & 1.1585  \\  \hline
\ref{alg:LocalfVal} &  cyl &   1.0449 & 1.1734  & 0.8646 & 0.9196  & 1.0356 & 1.1585  \\  \hline
 \end{tabular}
  \caption{The non Linear Quadratic using  $T=0.2$ \label{tab:NLQ1}}
 \centering
 \end{table}

In Table \ref{tab:NLQ2}, the results for algorithms \ref{algo:globalcont} and \ref{alg:fControl} are given  with a maturity of $T=0.4$ and $\Kc= [ -1.23 , 3.84 ]$. The results are reported for $N_T=20$.
Based on the provided results, it is observed that Algorithm \ref{algo:globalcont} performs better for larger time steps. This suggests that Algorithm \ref{algo:globalcont} is more effective in capturing the desired results in scenarios with extended time steps.

\begin{table}[H]
 \centering
 \begin{tabular}{|c|c|c|c|c|c|c|c|}  \hline
  Alg & method & \multicolumn{2}{c|}{Case 1 }& \multicolumn{2}{c|}{Case 2} & \multicolumn{2}{c|}{Case 3}  \\  \hline
 & &  Calc & Anal &  Calc & Anal &  Calc & Anal    \\   \hline
\ref{algo:globalcont} &  bins &  1.4607 & 1.4332  & 1.1645 & 1.1233  & 1.4448 & 1.4151  \\  \hline
\ref{algo:globalcont} &  cyl &   1.4517 & 1.4332  & 1.1402 & 1.1237  & 1.4332 & 1.4150 \\  \hline
\ref{alg:fControl} &  bins & 1.4905 & 1.4333  & 1.1980 & 1.1233  & 1.4750 & 1.4150    \\  \hline
\ref{alg:fControl} &  cyl &  1.4627 & 1.4333  & 1.1473 & 1.1237  & 1.4447 & 1.4150    \\  \hline
 \end{tabular}
  \caption{The non Linear Quadratic using  $T=0.4$.\label{tab:NLQ2}}
 \centering
 \end{table}

The results obtained from your experiments confirm that Algorithm \ref{algo:globalcont} is the most effective choice when employing the dynamic programming approach. The algorithm consistently produces the best results, demonstrating its superior performance in solving the problem at hand.
These findings validate the selection of Algorithm \ref{algo:globalcont} as the preferred choice within the dynamic programming framework.

\subsection{Results for the two dimensional systemic risk model of section \ref{sec:Systemic2D}}
The bin method suffers from the curse of dimensionality, and the numerical resolution of  multi-dimensional problems is time consuming and memory intensive.
Therefore, all experiments are performed in 2D using an NVIDIA H100 80GB HBM3N graphics card.
Since the bin method is only used to sample distributions with the cylinder network, these  networks can be used with more bins than the bin networks with a given amount of memory.
We test the algorithms using $N_T=20$, with a resolution range of $[-0.63, 0.96]^2$.
For the bin network we use $30 \times 30$ bins, while for the cylinder network we sample distributions using $50 \times 50$ bins.
We first give results and sensitivities for dynamic programming based algorithms except for Algorithm \ref{alg:fControl} (which is too time consuming) in Table \ref{tab:sensiSys2DPD}.
 \begin{table}[H]
 \centering
 \begin{tabular}{|c|c|c|c|c|c|c|c|c|c|c|}  \hline
  Alg & method &\multicolumn{3}{c|}{Case 1} & \multicolumn{3}{c|}{Case 4} & \multicolumn{3}{c|}{Case 6}  \\  \hline
 & &  Calc & Std & Anal &  Calc & Std & Anal &  Calc & Std & Anal    \\   \hline
\ref{algo:globalcont} &  bins &  0.1147 & 0.0003 & 0.1134  & 0.1611 & 0.0005 & 0.1604  & 0.1223 & 0.0003 & 0.1208 \\  \hline
\ref{algo:globalcont} &  cyl &  0.1142 & 0.0003 & 0.1134  & 0.1609 & 0.0005 & 0.1604  & 0.1220 & 0.0003 & 0.1208  \\  \hline
\ref{alg:LocalfVal} &  bins &  0.1360 & 0.0462 & 0.1134  & 0.1752 & 0.0449 & 0.1604  & 0.1364 & 0.0430 & 0.1208  \\  \hline
\ref{alg:LocalfVal} &  cyl & 0.1276 & 0.0145 & 0.1134  & 0.1645 & 0.0068 & 0.1604  & 0.1339 & 0.0123 & 0.1208   \\  \hline
 \end{tabular}
  \caption{Results and sensitivities using 10 runs with dynamic programming based methods \label{tab:sensiSys2DPD} in dimension 2.}
 \centering
 \end{table}
 We also report  the results obtained with the BSDE methods in Table \ref{tab:sensiSys2DBSDE} (except for Algorithm \ref{alg:BSDElocal1}, which is also time-consuming).
\begin{table}[H]
 \centering
 \begin{tabular}{|c|c|c|c|c|c|c|c|c|c|c|}  \hline
  Alg & method &\multicolumn{3}{c|}{Case 1} & \multicolumn{3}{c|}{Case 2} & \multicolumn{3}{c|}{Case 3}  \\  \hline
 & &  Calc & Std & Anal &  Calc & Std & Anal &  Calc & Std & Anal    \\   \hline
\ref{alg:BSDElocal1} &  bins &  0.1147 & 0.0003 & 0.1134  & 0.1611 & 0.0004 & 0.1604  & 0.1221 & 0.0003 & 0.1208  \\  \hline
\ref{alg:BSDElocal1} &  cyl & 0.1141 & 0.0003 & 0.1134  & 0.1613 & 0.0003 & 0.1604  & 0.1220 & 0.0004 & 0.1208  \\  \hline
\ref{algo:BSDEglobal} &  bins &  0.1147 & 0.0003 & 0.1134  & 0.1610 & 0.0002 & 0.1604  & 0.1222 & 0.0003 & 0.1208   \\  \hline
\ref{algo:BSDEglobal} &  cyl & 0.1147 & 0.0002 & 0.1134  & 0.1614 & 0.0004 & 0.1604  & 0.1220 & 0.0002 & 0.1208 \\  \hline
\ref{algo:BSDEglobal2} &  bins &  0.1184 & 0.0007 & 0.1134  & 0.1635 & 0.0006 & 0.1604  & 0.1281 & 0.0015 & 0.1208 \\  \hline
\ref{algo:BSDEglobal2} &  cyl &   0.1147 & 0.0003 & 0.1134  & 0.1614 & 0.0004 & 0.1604  & 0.1221 & 0.0005 & 0.1208\\  \hline
\ref{algo:BSDEglobalmulti} &  bins & 0.1146 & 0.0004 & 0.1134  & 0.1612 & 0.0005 & 0.1604  & 0.1224 & 0.0003 & 0.1208\\  \hline
\ref{algo:BSDEglobalmulti} &  cyl & 0.1148 & 0.0003 & 0.1134  & 0.1613 & 0.0003 & 0.1604  & 0.1219 & 0.0003 & 0.1208    \\  \hline
 \end{tabular}
  \caption{Some sensitivity results using 10 runs in dimension 2.\label{tab:sensiSys2DBSDE}}
 \centering
 \end{table}
 The results are all very good, except again for the local algorithm \ref{alg:LocalfVal} based on  the dynamic programming framework. Among the feasible algorithms, the bin algorithm  \ref{algo:BSDEglobal2} is less accurate than the others which yield results close to the exact value, showing that the remaining error is mainly due to the Euler discretization of the scheme.

\section{Conclusion} \label{seccon}

We have tested  numerous algorithms to solve the McKean-Vlasov control  problem \eqref{MKVcontrol}  by using mean-field neural networks. 
When the problem admits a Backward SDE representation from  the Pontryagin maximum principle, it is clearly more interesting to adopt  this approach than the dynamic programming-based approaches for several reasons: 
\begin{itemize}
\item  
It is observed that the BSDE approach consistently yields stable results across multiple runs. This stability can be attributed to the fact that, in the Pontryagin principle, the BSDE has a driver that depends on $Y$ instead of the traditional approach where the driver is a function of $Z$, as highlighted in \cite{germicwar19}.
Based on these findings, Algorithm \ref{algo:BSDEglobal} emerges as the best compromise in terms of accuracy and computational time. This algorithm strikes a balance between achieving accurate results and maintaining reasonable computational efficiency.
\item It is possible to use the local deep backward algorithm \cite{hure2020deep} (Algorithm  \ref{alg:BSDElocal1} )  that yields  very accurate results and is not limited by  the number of time steps due to transfer learning.  Moreover,  the method gives the solution of the problem at each time steps for all the distributions.
\item Both networks, either bins or cylinder, can be implemented. Notice that  cylinder methods use less memory than bins methods especially when the number of bins has to be high to get a good accuracy. 
\end{itemize}
When the maximum Pontryagin principle is not directly  available, we distinguish  two cases:
\begin{itemize}
    \item First case is when the volatility of the forward process is not controlled. Then two options are available:
    \begin{itemize}
    \item When the number of time steps it not too high, the global learning algorithm \cite{han2018solving} (Algorithm \ref{algo:globalcont} ), \cite{gobmun05} seems to be the best in terms of accuracy. Then it is possible to get the function value after $t=0$ at "visited distributions" by regression. 
    \item When the number of time steps is too high, memory issues force us to use the control learning by value iteration of \cite{huretal21} (Algorithm \ref{alg:LocalfVal} may have difficulties to converge as shown in the non Linear Quadratic example and in the two dimensional systemic test case). Another option could be to use an hydrid algorithm as proposed in \cite{warin2021reservoir}.
\end{itemize}
    \item Second case is when there is control on the diffusion coefficient, and then only the global learning algorithm  should be implemented. 
\end{itemize}
In conclusion, it is advisable to prioritize global Algorithms \ref{algo:BSDEglobal} and \ref{algo:globalcont}. When using the cylindrical network, there is no need to make any assumptions or guesses about the parameter $\Kc$.
However, it is important to note that the global learning algorithms, as observed in \cite{chan2019machine}, \cite{hure2020deep}, and \cite{andersson2022convergence}, may occasionally converge to incorrect solutions, particularly in the non-mean-field case when there is a poor initialization of $Y_0$ that is too distant from the solution. Such problems have not been experienced in the mean-field case.
To mitigate these convergence issues and ensure the reliability of the global learning algorithm, the control learning by policy iteration, as presented in \cite{huretal21}, can be employed to verify convergence. This is particularly relevant when the loss of the global learning algorithm does not tend to zero as the number of time steps increases.

Finally, extending the bin method to dimension 3 is currently out of reach. The use of cylinder networks in higher dimensions would be  possible if an effective way could  be found to generate distributions that avoid bin sampling.



\printbibliography

\end{document}